\documentclass[11 pt,a4paper,twoside]{article}
\usepackage{amsmath,amssymb,amsthm,amsfonts,epsfig,graphics}

\renewcommand{\AA}{\mathbb{A}}

\newcommand{\DD}{\mathbb{D}}

\newcommand{\NN}{\mathbb{N}}

\newcommand{\RR}{\mathbb{R}}
\renewcommand{\SS}{\mathbb{S}}
\newcommand{\TT}{\mathbb{T}}

\newcommand{\ZZ}{\mathbb{Z}}

\def\cD{{\cal D}}

\renewcommand{\phi}{\varphi}
\renewcommand{\epsilon}{\varepsilon}

\newtheorem{theo}{Theorem}[section]
\newtheorem{prop}[theo]{Proposition}
\newtheorem{coro}[theo]{Corollary}
\newtheorem*{conj*}{Conjecture}
\newtheorem{lemma}[theo]{Lemma}
\newtheorem{sublemma}[theo]{Sublemma}
\newtheorem{fact}[theo]{Fact}

\newtheorem{rema}[theo]{Remark}

\renewcommand{\t}{\widetilde}
\def\wt{\widetilde}
\def\tA{\wt \AA}

\newtheorem*{theo*}{Theorem}

\def\inte{\mathrm{Int}}
\def\rot{\mathrm{Rot}}

\newcommand{\dd}{\operatorname{d}}
\newcommand{\homeo}{\operatorname{Homeo}}
\newcommand{\id}{\operatorname{Id}}
\newtheorem{lemmanota}[theo]{Lemma and notation}

\title{Pseudo-rotations of the closed annulus: \\variation on a 
theorem of J.\,Kwapisz}
\author{F. B\'eguin\footnote{Laboratoire de math\'ematiques, Universit\'e Paris Sud,
91405 Orsay Cedex, France.},
 S. Crovisier\footnote{Laboratoire de Topologie, Universit\'e de Bourgogne, BP
47 870, 21078 Dijon Cedex, France.},
 F. Le Roux\footnote{Laboratoire de math\'ematiques, Universit\'e Paris Sud,
91405 Orsay Cedex, France.}
 and A. Patou\footnote{Laboratoire de Math\'ematiques, Universit\'e
Grenoble 1,
 BP 74, 38402 Saint-Martin-d'H\`eres Cedex, France.}
}

\begin{document}

\sloppy

\maketitle
% abstract et intro
%abstract coquille ``of angle...''
%-- pb d'isotope a identit\'e

%-- $f$ / $h$

%-- ; :
 
\begin{abstract}
Consider a homeomorphism $h$ of the closed annulus $\SS^1\times 
[0,1]$, isotopic to the identity, such that the rotation
set of $h$ is reduced to a single irrational number $\alpha$ (we say that $h$ 
is an \emph{irrational pseudo-rotation}).
 For every positive integer $n$, we prove that there exists a 
simple arc $\gamma$ joining one of the boundary component of
the annulus to the other one, such that $\gamma$ is disjoint from its $n$ 
first iterates under $h$.  As a corollary, we obtain that
 the rigid rotation of angle $\alpha$ can be approximated by
homeomorphisms conjugate to $h$.
%there exists homeomorphisms conjugate to $h$ arbitrarily close
%to the rigid rotation of angle $\alpha$.
% $h$ can be conjugate to a homeomorphism arbitrarily close to a rigid 
%rotation. 
The first result stated above is an analog of a
theorem of J.\,Kwapisz dealing with 
%homeomorphisms
diffeomorphisms of the two-torus; 
we give some new, purely two-dimensional, proofs,
that work both for the annulus and for the torus case.


\end{abstract}

\section{Introduction}

The concept of \emph{rotation number} was introduced by H. Poincar\'e 
to study the dynamics of circle
homeomorphisms (in the context of torus flows, see \cite{Poi},
chapitre XV).
 More precisely, for every orientation-preserving 
homeomorphism~$h$ of the circle
$\SS^1=\RR/\ZZ$, Poincar\'e defined an element of $\SS^1$,  measuring 
the ``asymptotic speed of rotation of  the orbits of
$h$ around the circle":  the so-called \emph{rotation number of $h$}. 
The central question in  this theory is: how much does
the dynamics of an orientation-preserving homeomorphism of the circle 
of rotation number $\alpha$ look like the dynamics
of a rigid rotation $R_\alpha$? The classical results obtained by 
Poincar\'e and A.\,Denjoy (\cite{Den}) provide a quite comprehensive
list of answers to this question:

\smallskip

\noindent $\bullet$ If $\alpha=p/q$ (where $p,q$ are relatively prime 
integers), then $h$ has at least one
periodic orbit, all the periodic orbits of $h$ have prime period $q$, 
and the cyclic order of the points of any periodic orbit of $h$
is the same as the cyclic order of the points of an orbit of the 
rotation $R_\alpha$. If $\alpha$ is irrational, then $h$
does not have any periodic orbit, and the cyclic order of the points 
of any orbit of $h$ is the same as the cyclic order of the
points of an orbit of the rotation $R_\alpha$.

\smallskip

\noindent $\bullet$ If $\alpha$ is irrational, then $h$ is 
semi-conjugate to the rotation $R_\alpha$;
moreover, $h$ is in the closure of the conjugacy class of the 
rotation $R_\alpha$, and $R_\alpha$
is in the closure of the conjugacy class of $h$ (\emph{i.e.} $h$ can 
be conjugated to a homeomorphism arbitrarily close to
the rotation $R_\alpha$, and the rotation $R_\alpha$ can be conjugated 
to a homeomorphism arbitrarily close to $h$).

\smallskip

\noindent $\bullet$ If $\alpha$ is irrational and $h$ is a 
$C^2$-diffeomorphism, then $h$ is conjugate to the rotation
$R_\alpha$.

\bigskip

Poincar\'e's construction of the rotation number can be generalized 
for homeomorphisms of the closed
annulus $\AA:=\SS^1\times [0,1]$. 
%The main difference is that, 
%instead of a single number, we obtain an interval. 
More precisely, for every homeomorphism $h$
of the closed annulus $\AA$ which is isotopic to the identity, the 
\emph{rotation set of $h$} is a closed interval of $\RR$,
\footnote{Contrary to what happens in the case of the circle, 
two different orbits of a homeomorphism of
the annulus might have different ``asymptotic speeds of rotation''.}
defined up to the addition of an integer, which measures the 
asymptotic speeds of rotation of the orbits of $h$
around the annulus (see section \ref{ss.rotation}). In the present article, 
we focus on homeomorphisms whose rotation set is a
``small'' interval. In particular, we call \emph{irrational 
pseudo-rotation} every homeomorphism of the closed annulus
$\AA$, isotopic to the identity, whose rotation set is reduced to a 
single irrational number $\alpha$ (and we say that
$\alpha$ is the \emph{angle} of the pseudo-rotation). In this 
context, the natural question is: how much does the dynamics of
an irrational pseudo-rotation of angle $\alpha$ look like the 
dynamics of the rigid rotation of angle $\alpha$? The aim of
the present article is to give some partial answer to this question.

\bigskip

We define an \emph{essential simple arc} in the annulus $\AA$ as a 
simple arc in $\AA$ joining one of the boundary
components of $\AA$ to the other one. We shall prove the following 
theorem (which is a variation on a result of J. Kwapisz,
dealing with torus diffeomorphisms, see \cite{Kwa1}):

\begin{theo}
\label{th.pseudo-rotation}
Let $h:\AA\rightarrow\AA$ be an irrational pseudo-rotation of angle $\alpha$.
Then, for every  positive integer $n$, there exists an essential 
simple arc $\gamma_n$ in $\AA$, such that the
arcs $\gamma_n,\dots,h^n(\gamma_n)$ are pairwise disjoint.
Moreover, the cyclic order of these arcs is the same as the cyclic 
order of the $n$ first iterates of a vertical segment $\{\theta\} \times[0,1]$
under the rigid rotation of angle $\alpha$.
\end{theo}

Actually, theorem \ref{th.pseudo-rotation} will appear as a corollary 
of a more technical statement, concerning
the homeomorphisms of the annulus whose rotation set is a ``small" 
interval (more precisely, a \emph{Farey interval},
see theorem \ref{th.closed-annulus}).

Theorem \ref{th.pseudo-rotation} can be considered as a 
two-dimensional version of the above mentioned result concerning the
cyclic order for circle homeomorphisms.
%fact that, for any circle
%homeomorphism of irrational rotation number $\alpha$, the cyclic 
%order of the points of any orbit is the same as the
%cylic order of the points of an orbit of the rigid rotation of angle 
%$\alpha$.
 Nevertheless, the situation is quite more
complicated than in the circle. For example, it is should be noticed 
that the statement of theorem \ref{th.pseudo-rotation} is
optimal  in the sense that it is impossible to make the arc 
$\gamma_n$ independent of $n$. More precisely, M. Herman
has constructed a $C^\infty$ irrational pseudo-rotation $h$ of the closed 
annulus $\AA$ which is not conjugate to a rigid rotation~; 
 it is not difficult to see that no essential simple  arc is
disjoint of all its iterates under $h$ (see \cite{Her} and \cite{Han}).

As we have already said, theorem \ref{th.pseudo-rotation} is a variation of an analog 
result of Kwapisz, dealing with diffeomorphisms
of the torus $\TT^2$. The true aim of our article is not to adapt 
Kwapisz's proof to the case of annulus
homeomorphisms, but rather to provide some completely different and 
(in our opinion) more natural proofs. Indeed, in his
proof, Kwapisz introduces the suspension of the diffeomorphism under 
consideration, and uses some
$3$-dimensional topology techniques to find the wanted curve as the 
intersection of two
cross-sections of this suspension. The two proofs of theorem \ref{th.pseudo-rotation} 
that we give in the present article are purely
two-dimensional, and only involve some classical manipulations on 
arcs. By the way, these proofs also work in the
torus case (see appendix \ref{app.torus}).

\bigskip

As a corollary of theorem \ref{th.pseudo-rotation}, we obtain the 
following result\footnote{It might be
interesting to notice that corollary \ref{c.closure-conjugacy-class} 
is actually equivalent to theorem
\ref{th.pseudo-rotation}.}:
\begin{coro}
\label{c.closure-conjugacy-class}
Let $h:\AA\rightarrow\AA$ be an irrational pseudo-rotation of angle 
$\alpha$. Then the rigid rotation $R_\alpha$ of angle
$\alpha$ is in the closure of the conjugacy class of $h$.
\end{coro}

In other words, any irrational pseudo-rotation can be conjugated to 
obtain a homeomorphism which is
arbitrarily close to a rigid rotation. Nevertheless, we point out 
that
% under the hypothesis of corollary \ref{c.closure-conjugacy-class},
 we do not know if the converse is 
true, namely,
if any irrational pseudo-rotation of angle $\alpha$ 
is in the closure of the conjugacy class of the
rotation $R_\alpha$.

Corollary \ref{c.closure-conjugacy-class} was motivated by the 
situation on the two-torus. Indeed, an analog of
corollary \ref{c.closure-conjugacy-class} holds on the two-torus; it 
is actually an immediate consequence of another
theorem by Kwapisz, called the \emph{tiling theorem} (see 
\cite{Kwa2}). The tiling theorem asserts roughly that if the
rotation set of a two-torus diffeomorphism $h$ is reduced to a single 
irrational point, then for any $n$, one can find a finite
tiling, which is almost invariant under $h$ (there are only three 
tiles whose images do not fit in with the tiling), and such that
the restriction of $h$ to the 1-skeleton is conjugate to the 
restriction of the corresponding rigid rotation to the
1-skeleton of a similar tiling. In the case of the annulus, theorem 
\ref{th.pseudo-rotation} also provides a kind of
almost invariant tiling of the annulus. Nevertheless, corollary 
\ref{c.closure-conjugacy-class} is a little more
difficult to derive in the annulus case because, unlike what happens in the 
torus setting, the diameter of the tiles of
the corresponding tiling for the rigid rotation does not go to zero 
when the number of tiles increase.

\bigskip

%Generalizations of Poincar\'e-Birkhoff theorem obtained by J.\,Franks 
%(\cite{Fra}) or C.\,Bonatti and L.\,Guillou
%(\cite{Gui}) imply that an homemorphism of the closed annulus, 
%isotopic to the identity, preserving the
%Lebesgue measure, is an irrational pseudo-rotation if and only if it 
%does not have any periodic orbit. The results
%of Franks and Bonatti-Guillou are actually more general ; putting 
%these results together with theorem
%\ref{th.pseudo-rotation}, we obtain the following noteworthy corollary:
Finally, it is interesting to associate
theorem~\ref{th.pseudo-rotation} with some generalizations of
the Poincar\'e-Birkhoff theorem obtained by J.\,Franks  
(\cite{Fra}) or C.\,Bonatti and L.\,Guillou (\cite{Gui}). Let us
recall the result of Bonatti and Guillou. It deals with homeomorphisms
$h$ of the closed  annulus $\AA$ 
that are isotopic to the identity, and claims that
 \emph{if $h$ is fixed point free,
then either  there exists an essential
simple arc in $\AA$ that is disjoint from its image under $h$,
or there exists a non-homotopically trivial 
simple closed curve in $\AA$ which is
disjoint from its image under $h$}. In particular, it implies that if
$h$ preserves the Lebesgue measure and has no periodic point, then $h$
is an irrational pseudo-rotation. Putting 
the quoted result together with theorem
\ref{th.pseudo-rotation}, we obtain the following corollary:
\begin{coro}
\label{c.without-periodic-orbit}
Let $h$ be a homeomorphism of the annulus $\AA$, which is isotopic 
to the identity, and which does not have
any periodic point. Then:

\smallskip

\noindent \emph{(i)} either there exists a non-homotopically trivial 
simple closed curve in $\AA$, which is
disjoint from its image under $h$,

\smallskip

\noindent \emph{(ii)} or  $h$ is an irrational 
pseudo-rotation, and, for every positive integer $n$, there
exists an essential simple arc $\gamma_n$ in $\AA$, such that the 
arcs $\gamma_n,h(\gamma_n),\dots,h^n(\gamma_n)$
are pairwise disjoint.
\end{coro}

%Observe that case (i) of corollary \ref{c.without-periodic-orbit} 
%cannot occur if the homeomorphism $h$
%preserves the measure the Lebesgue measure of the annulus $\AA$.

\bigskip

In a forthcoming paper, we shall prove some analogs of theorem 
\ref{th.pseudo-rotation}, corollaries
\ref{c.closure-conjugacy-class} and \ref{c.without-periodic-orbit} 
for homeomorphisms of the open annulus
$\SS^1\times ]0,1 [$. Some completely different (and more 
sophisticated) proofs are needed. All the difficulty arise
from the lack of compacity of the open annulus, which forces to 
change the definition of the rotation set (in particular,
one has to restrict to measure-preserving homeomorphisms).

{ \bf Acknowledgements.} We would like to thank Jarek Kwapisz for helpfull discussions on his work.
\section{Definitions and precise statement}

\subsection{Notations}

In this paper, we denote by $\AA$ the closed annulus $\SS^1\times 
[0,1]$, and  by
$\wt\AA:=\RR\times [0,1]$ the universal covering of $\AA$. Moreover, 
we denote by $\pi$ the
canonical projection of $\wt\AA$ onto $\AA$, and  by 
$T:\wt\AA\rightarrow\wt\AA$
the translation defined by $T : (\theta,t) \rightarrow (\theta+1,t)$.

Of course, the translation $T$ is a generator of the automorphism 
group of the projection
$\pi:\wt\AA\rightarrow\AA$. We observe that if $h$ is a 
homeomorphism of the annulus
$\AA$ which is isotopic to the identity, and if $\wt h$ is a lift of 
$h$ to the band $\tA$, then
$\wt h$ commutes with the covering translation $T$. Conversely, every 
homeomorphism of
$\tA$ which is isotopic to the identity and which 
 commutes with $T$ is the lift of  a homeomorphism of the 
annulus $\AA$ isotopic
to the identity.
For every $\alpha\in\RR$, we denote by $R_\alpha: (\theta,t)  \mapsto
(\theta+\alpha,t) $ the rigid rotation of angle $\alpha$ in the
annulus $\AA$.
%, and we denote by $T_\alpha$ the translation of vector 
%$(\alpha,0)$ in the
%band $\wt\AA$, that is
%$$\begin{array}{llllllll}
%R_\alpha : & \AA & \rightarrow & \AA & \quad \quad T_\alpha: & \AA & 
%\rightarrow & \AA\\
%   & (\theta,t) & \mapsto & (\theta+\alpha,t) & & (\theta,t) & \mapsto 
%& (\theta+\alpha,t)
%\end{array}$$
Finally, we denote by $p_1 : \tA=\RR \times [0,1] \rightarrow \RR$ be the first
coordinate projection.

\subsection{The rotation set of a homeomorphism of the annulus $\AA$}
\label{ss.rotation}
Let $h$ be a homeomorphism of the bounded annulus $\AA$ which is
isotopic to the identity, and $\wt h:\wt\AA\rightarrow\wt\AA$ be a lift of $h$.
We define the $n^{th}$ \emph{displacement set} of $\wt h$ to be the set
$$
D_n(\wt h)=\left\{\frac{p_1(\wt h^n (\wt x)) -p_1(\wt x)}{n} \mid \wt 
x \in \tA \right\}.
$$
Then a real number $v$ is called a \emph{rotation vector} if it is the
limit of a sequence $(v_{n_k})_{k \geq 0}$ such that each $v_{n_k}$ belongs to
the $n_k^{th}$ displacement set of $\wt h$, where the sequence $(n_k)$ tends
to $+\infty$. The \emph{rotation set of $\wt h$} is the set $\rot(\wt 
h)$ of all rotation
vectors. It is easy to see that the sets $D_n(\wt h)$ and $\rot(\wt 
h)$ are compact intervals.
This definition of the rotation set of $\wt h$ is the analog of a 
definition given by
Misiurewicz and Zieman in \cite{MisZie} in the case of torus homeomorphisms.

Let us recall briefly an alternative definition that follows an idea of
S. Schwartzman \cite{Sch}. If $\mu$ is an
invariant measure for $\wt h$, the \emph{rotation vector of $\mu$} is
$$\int_D \left( p_1(\wt h (\wt x)) -p_1(\wt x) \right) d \mu (\wt x)$$
where $D$ is any fundamental domain of the covering $\tA$ (for example
 $D=[0,1[ \times [0,1]$). If $\mu$ is ergodic, then the ergodic
theorem implies that the rotation number $v$ of $\mu$ is realized, in the
following strong sense: there exists a point $\wt x$ such that
$$\lim_{n \rightarrow +\infty} \frac{p_1(\wt h^n (\wt x)) -p_1(\wt x)}{n}= v.$$
One can deduce from this that the rotation set of $\wt h$ coincides
with the set of rotation vectors of all invariant measures, and
  that the endpoints of the interval $\rot(\wt h)$ are realized in the 
sense defined above.

For any integers $p,q$, the map $\wt h^q \circ
T^{-p}$ is a lift of $h^q$. Using the fact that $\wt h$ and $T$ 
commute, we have the following
easy property:

\begin{lemma}
\label{lemm.rotation-itere}
For any couple $(p,q)$ of integers, the rotation set of $\wt h^q \circ
T^{-p}$ is given by
$$\rot(\wt h^q \circ T^{-p})= q \times \rot(\wt h) -p.$$
\end{lemma}

In particular, the rotation sets of two different lifts differ by an
integer, so that the rotation set of $h$ is well defined as an
interval of $\RR$ modulo $\ZZ$ (formally, we can see it as an element
of $\RR^2$ quotiented by the action of $(v_1,v_2) \mapsto (v_1+1,v_2+1)$).
It is an invariant with respect to the conjugacy by the homeomorphisms
of $\AA$ that are isotopic to the identity.

%%%%%%
\subsection{Cyclic order on the circle and the annulus}
The natural orientation of $\RR$ induces an orientation on the circle
$\SS^1=\RR/\ZZ$. Given three distinct points $p_1,p_2,p_3$ on $\SS^1$,
we will say that the triplet $(p_1,p_2,p_3)$ is \emph{positive} if
the point $p_2$ is crossed when going from $p_1$ to $p_3$ in the
positive way.

A simple arc $\gamma : [0,1] \rightarrow \AA$ is called an \emph{essential
simple arc} if $\gamma$ joins one of the boundary components of $\AA$ 
to the other, and if
$\gamma(]0,1[)$ is included in the interior of $\AA$. We define 
similarly the notion of
\emph{essential simple arc} in the band $\wt\AA$. Similarly to the 
cyclic order for distinct
points on the circle, we define a cyclic order on triplets of essential 
simple arcs in $\AA$ that are pairwise disjoint. Note that this
can be done simply by considering the endpoints of the three arcs on
one of the boundary components. We will also use the (related) total order on
sets of  pairwise disjoint essential 
simple arcs in $\wt \AA$.

%%%%%%%
\subsection{Farey intervals}

In this article, all the rational numbers $\frac{p}{q}$ (with 
$p\in\ZZ$ and $q\in\NN\setminus\{0\}$)
will be written as an irreducible fraction. A \emph{Farey interval} 
is an interval
$]\frac{p}{q},\frac{p'}{q'}[$ of
$\RR$ with rational endpoints, such that $p'q-pq'=1$ (which amounts to
saying that the length of the interval is $\frac{1}{qq'}$). Some
elementary properties of Farey intervals are used in sections \ref{s.first-proof}
 and \ref{s.alternative-proof} and proved in 
appendix~\ref{app.arithmetic}.

%%%%%%%%
\subsection{Precise statement of the theorem}

\begin{theo}
\label{th.closed-annulus}
Let $h:\AA\rightarrow\AA$ be a homeomorphism isotopic to the identity, and
let $\widetilde h:\widetilde\AA\rightarrow\widetilde\AA$ be a lift of
$h$. Assume that the rotation set of
$\widetilde h$ is included  in a Farey interval
$]\,\frac{p}{q},\frac{p'}{q'}\,[$. Then there exists an essential
simple arc $\gamma$ in
the annulus $\AA$ such that the arcs
$\gamma,h(\gamma),\dots,h^{q+q'-1}(\gamma)$ are pairwise disjoint.
Moreover the cyclic order of these arcs in $\AA$ is the same as the 
cyclic order of the first $n$ iterates of a vertical
segment under any rigid rotation of angle $\alpha\in 
]\,\frac{p}{q},\frac{p'}{q'}\,[$.
\end{theo}

The statement on cyclic order means that
 for any $k_1,k_2,k_3 \in \{0, \dots, q+q'-1\}$,
 one has the equivalence
$$
(h^{k_1}(\gamma),h^{k_2}(\gamma),h^{k_3}(\gamma))
\mbox{ is positive } \Leftrightarrow (k_1 \alpha, k_2 \alpha, k_3
\alpha) \mbox{ is positive}
$$
where the numbers $k_i \alpha$ are considered as elements of the circle $\RR/\ZZ$.

The result announced in the introduction (theorem
 \ref{th.pseudo-rotation}) follows 
directly from this new one, by
noting that given any irrational number $\alpha$, one can find a Farey interval
$]\,\frac{p}{q},\frac{p'}{q'}\,[$ containing $\alpha$ with $q+q'$ 
arbitrarily big.

%%%%%%
\subsection{A basic property}

\begin{lemma}
\label{lemm.ca-pousse}
Suppose  that the rotation set of $\wt h$ is included in
$]0,+\infty[$, and choose a real number $\rho$ such that $0
< \rho < \inf(\rot(\wt h))$.
Then there exists a real number $s$ such that for every $\wt x$ in 
$\tA$, for every
positive integer $n$,
\begin{equation}\label{e.fact}
p_1(\wt h^n(\wt x)) \geq p_1(\wt x) +\rho n -s.
\end{equation}
\end{lemma}

The proof is easy, and left to the reader. We shall use the following 
consequence of
this lemma:  \emph{under the hypotheses of lemma \ref{lemm.ca-pousse}, for every
compact subset $K$ of $\t A$  and every $s_0>0$, there exists $n_0>0$ 
such that for all
$n \geq n_0$, $\wt h^n(K) \subset [s_0, +\infty[ \times [0,1]$}.

\begin{rema}
Another consequence is that, under the hypotheses of lemma 
\ref{lemm.ca-pousse}, the quotient
space $\wt\AA/\wt h$ is separated. Using the classification of 
surfaces, one can see that
$\wt\AA/\wt h$ is necessarily homeomorphic to a closed annulus, which 
implies that $\wt h$ is
conjugate to a translation. We shall not need this fact.
\end{rema}
%
% ajouter les defs d'arcs simple essentiels en haut et en bas
%parentheses/crochets

\section{First proof of the main theorem}
\label{s.first-proof}
%In this section, we give our first proof of theorem
%\ref{th.closed-annulus}. 
Two proofs of  theorem \ref{th.closed-annulus} will be given,
the first one in this section and the second one in the following
section. These sections can be read in any order.

The first proof uses two kinds of ingredients:
some elementary arithmetical properties of Farey intervals, and some 
(classical) operations on
essential simple arcs in the band $\widetilde\AA$.

%%%%%%%%%%%%%%%%%%%%%%%
\subsection{Some more notations}

For every essential simple arc $\Gamma$ in the band $\widetilde\AA$, 
we denote by $R(\Gamma)$ the
closure of the connected component of $\widetilde\AA\setminus\Gamma$ 
which is ``on the right" of the arc
$\Gamma$.

Given an essential simple arc $\Gamma$ in $\widetilde\AA$ and a homeomorphism
$\Psi:\widetilde\AA\rightarrow\widetilde\AA$,  we say that the set 
$R(\Gamma)$ is an \emph{attractor}
(resp. a \emph{strict attractor}) for $\Psi$ if the image of 
$R(\Gamma)$ under $\Psi$ is included in
$R(\Gamma)$ (resp. in the interior of $R(\Gamma)$). Observe that if 
$\Psi$ is isotopic to the identity, we have
$R(\Psi(\Gamma))=\Psi(R(\Gamma))$, so that the set
$R(\Gamma)$ is a (strict) attractor for $\Psi$ if and only if the 
image of the arc $\Gamma$ under $\Psi$ is
included in (the interior of) $R(\Gamma)$.

%%%%%%%%%%%%%%%%%
\subsection{Attractors for families of commuting homeomorphisms}

Theorem \ref{th.closed-annulus} will follow from elementary arithmetical 
properties of Farey intervals, and from the
following proposition:
\begin{prop}
\label{p.plusieurs-applications}
Let $\Psi_0,\dots,\Psi_p$ be some homeomorphisms of the  band
$\widetilde\AA$ isotopic to the identity, 
pairwise commuting, and commuting with the translation $T$. Assume
that the rotation set of each of these homeomorphisms 
 is included in $]0,+\infty[$. Then there exists an essential simple arc $\Gamma$ in 
$\widetilde\AA$,
such that the set $R(\Gamma)$ is a strict attractor for each of the 
homeomorphisms
$\Psi_1,\dots,\Psi_p$.
\end{prop}

%The next three pages are devoted to the proof of proposition 
%\ref{p.plusieurs-applications}. 
We begin with a technical point which consists in  describing 
an operation on essential simple arcs. This operation will be used
intensively to construct the simple arc demanded by
proposition~\ref{p.plusieurs-applications}. The proof of the following lemma is
postponed to section~\ref{s.lvee}.

\begin{figure}[ht]
\par \centerline{\hbox{\input{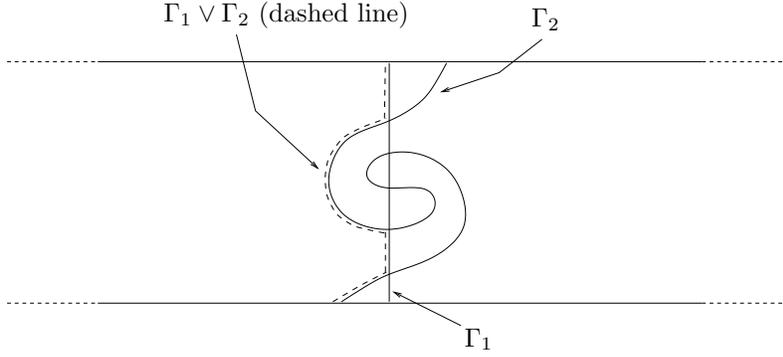}}} \par
\caption{\label{f.vee}Two essential simple arcs $\Gamma_1$ and 
$\Gamma_2$, and the arc
$\Gamma_1\vee\Gamma_2$}
\end{figure}

\begin{lemmanota}[figure~\ref{f.vee}]
\label{l.vee}
Let $\Gamma_1$ and $\Gamma_2$ be two essential simple arcs in 
$\widetilde\AA$, and let $U$ be
the  unique non-bounded connected component of the set 
$(\widetilde\AA\setminus R(\Gamma_1))\cap
(\widetilde\AA\setminus R(\Gamma_2))$. The boundary  (in 
$\widetilde\AA$) of $U$ is an essential simple
arc, that we denote by $\Gamma_1\vee\Gamma_2$.
\end{lemmanota}

\begin{rema}
\label{r.vee}
Let $\Gamma_1$ and $\Gamma_2$ be two essential simple arcs in 
$\widetilde\AA$. The following
properties are immediate consequences of the definition of the arc 
$\Gamma_1\vee\Gamma_2$:

\smallskip

\noindent \textbf{\emph{(i)}} The arc $\Gamma_1\vee\Gamma_2$ is 
included in the union of the arcs
$\Gamma_1$ and $\Gamma_2$.

\smallskip

\noindent \textbf{\emph{(ii)}} The sets $R(\Gamma_1)$ and 
$R(\Gamma_2)$ are included in the set
$R(\Gamma_1\vee\Gamma_2)$.
\end{rema}

\begin{rema}
The operation which maps two essential simple arcs 
$\Gamma_1,\,\Gamma_2$ to the essential
simple arc $\Gamma_1\vee\Gamma_2$ is associative (and commutative). 
Therefore, given any finite
number of essential simple arcs $\Gamma_1,\dots,\Gamma_n$, the essential simple
arc $\Gamma_1\vee\dots\vee\Gamma_n$ is well-defined.
\end{rema}

Now we are in a position to prove  proposition
\ref{p.plusieurs-applications}.

\begin{proof}[Proof of proposition \ref{p.plusieurs-applications}]
We proceed by induction. For every $k\in\{1,\dots,p\}$, we shall 
construct an essential simple arc
$\Gamma_k$ such that the set $R(\Gamma_k)$  is a strict attractor for 
the homeomorphisms
$\Psi_1,\dots,\Psi_k$.

\paragraph{Construction of the arc $\Gamma_1$.} Let $\Gamma_0$ be an 
essential simple arc in
$\widetilde\AA$. According to lemma \ref{lemm.ca-pousse}, there exists an integer 
$N$ such that  the
arc $\Psi_1^N(\Gamma_0)$ is included in $R(\Gamma_0)$. 
%(or equivalently, such that the set
%$R(\Psi_1^N(\Gamma_0))$ is included in $R(\Gamma_0)$). 
We consider the essential simple arc
$$\Gamma_1:=\Gamma_0 \vee\Psi_1(\Gamma_0) \vee \dots\vee\Psi_1^{N-1}(\Gamma_0).$$
By item (ii) of remark \ref{r.vee}, the sets 
$R(\Gamma_0),\dots,R(\Psi_1^{N-1}(\Gamma_0))$  are included
in $R(\Gamma_1)$. Moreover, by definition of the integer $N$, the set 
$R(\Psi_1^N(\Gamma_0))$ is
included in $R(\Gamma_0)$, and therefore it is also included in
$R(\Gamma_1)$.  In particular, the arcs 
$\Gamma_0,\dots,\Psi_1^N(\Gamma_0)$ are
included in $R(\Gamma_1)$. On the other hand, by item (i) of remark 
\ref{r.vee}, the arc $\Gamma_1$ is
included in the union of the arcs 
$\Gamma_0,\dots,\Psi_1^{N-1}(\Gamma_0)$~; therefore the arc
$\Psi_1(\Gamma_1)$ is included in the union of the arcs 
$\Psi_1(\Gamma_0),\dots,\Psi_1^{N}(\Gamma_0)$.
Putting everything together, we obtain that the arc 
$\Psi_1(\Gamma_1)$ is included in the set
$R(\Gamma_1)$. Hence the set
$R(\Gamma_1)$ is an attractor for the homeomorphism $\Psi_1$.
It remains to perturb $\Gamma_1$  in such a way that the set $R(\Gamma_1)$ 
becomes a strict  attractor for $\Psi_1$~; this is made possible by
the following  technical lemma~:
\begin{lemma}
\label{l.perturbation}
Let $\Psi:\widetilde\AA\rightarrow\widetilde\AA$ be a
homeomorphism  commuting with the
translation $T$, such that the rotation set of $\Psi$ is included in
$]0,+\infty[$. Suppose that we have found an essential
simple arc $\Gamma$ in $\widetilde \AA$, such that $R(\Gamma)$ is an
attractor for $\Psi$. Then, arbitrarily close to
$\Gamma$, there exists an essential simple arc $\widehat\Gamma$ such
that $R(\widehat\Gamma)$ is a strict attractor
for $\Psi$.
\end{lemma}
This lemma is extracted from (\cite[part 5]{Gui})~; a (slight)
variation on the proof of \cite{Gui} is given in section~\ref{ss.perturbation}.

\paragraph{Induction step.} Let $k\in\{1,\dots,p-1\}$. We assume that 
we have constructed an essential
simple arc $\Gamma_k$ such that the set $R(\Gamma_k)$ is a strict 
attractor for the homeomorphisms
$\Psi_1,\dots,\Psi_k$. By lemma \ref{lemm.ca-pousse}, there exists an integer $N$ 
such that the
arc $\Psi_{k+1}^N(\Gamma_k)$ is included in $R(\Gamma_k)$. We 
consider the essential simple
arc
$$\Gamma_{k+1}:=\Gamma_k\vee \Psi_{k+1}(\Gamma_k)\vee\dots\vee \Psi_{k+1}^{N-1}(\Gamma_k).$$
The same argument as in the construction of the arc $\Gamma_1$ shows 
that the set $R(\Gamma_{k+1})$ is
an attractor for  the  homeomorphism $\Psi_{k+1}$. Now, let $j\in 
\{1,\dots,k\}$; we will check that the set $R(\Gamma_{k+1})$ is still
a strict attractor for $\Psi_j$ (this will essentially follow from the
fact that $\Psi_j$ commutes with $\Psi_k$). Firstly, by item (i) of remark
\ref{r.vee}, the arc $\Gamma_{k+1}$ is included
in the union of the arcs $\Gamma_k,\dots,\Psi_{k+1}^{N-1}(\Gamma_k)$. 
%Hence,  the arc
%$\Psi_j(\Gamma_{k+1})$ is included in the union of the arcs
%$\Psi_j(\Gamma_k),\dots,\Psi_j(\Psi_{k+1}^{N-1}(\Gamma_k))$. 
Secondly,  we 
observe that the sets $R(\Gamma_k),\dots,R(\Psi_{k+1}^{N-1}(\Gamma_k))$ are
strict attractors for the homeomorphism $\Psi_j$ (this is because the set 
$R(\Gamma_k)$ is a strict attractor for
the homeomorphism $\Psi_j$, and because the homeomorphisms $\Psi_k$ 
and $\Psi_j$ commute).
Hence, for every $j\in \{1,\dots,k\}$, we have
$$
\Psi_j(\Gamma_{k+1}) \subset \bigcup_{i=0}^{N-1} 
\Psi_j(\Psi_{k+1}^i(\Gamma_k))\subset
\bigcup_{i=0}^{N-1} \inte(R(\Psi_{k+1}^i(\Gamma_k)))\subset 
\inte(R(\Gamma_{k+1}))
$$
where the last inclusion comes from both items of remark~\ref{r.vee}.
So the set $R(\Gamma_{k+1})$ is a strict attractor for the 
homeomorphisms $\Psi_1,\dots,\Psi_k$.
Then, using lemma \ref{l.perturbation}, we can perturb the arc 
$\Gamma_{k+1}$ in such a way that the
set $R(\Gamma_{k+1})$ becomes a strict attractor for the 
homeomorphism $\Psi_{k+1}$. Provided
that the perturbation is small enough, we keep the property that 
$R(\Gamma_{k+1})$ is
a strict attractor for the homeomorphisms $\Psi_1,\dots,\Psi_k$ (being a 
strict attractor is an ``open property"). This completes the proof of
proposition~\ref{p.plusieurs-applications}.
\end{proof}

%%%%%%%%%%%%%%%%%%%%%%
\subsection{Proof of the theorem}
We now turn to the proof of the main theorem.
It consists in applying proposition 
\ref{p.plusieurs-applications} to a well-chosen family of
homeomorphisms $\Psi_1,\dots,\Psi_{q+q'}$. Each of these homeomorphisms
 will be obtained as the composition of a power of $\wt h$ and a power of $T$.

\begin{proof}[Proof of theorem \ref{th.closed-annulus}]
 Let $\rho$ be any number in the Farey interval 
$]\,\frac{p}{q}\,,\,\frac{p'}{q'}\,[$. According to
section~\ref{ssapp.arithmetic} of the appendix, we have:
\begin{itemize}
\item For every $k\in \{1,\dots,q+q'-1\}$, the number $k.\rho$ is not
an integer, so we may  define the number $\alpha_k \in ]0,1[$ and the
integer $n_k$ such that $k.\rho=n_k+\alpha_k$;
\item the numbers $\alpha_1, \dots, \alpha_{q+q'-1}$ are distinct, so 
we may consider the permutation $\sigma$ of the set $\{1,\dots,q+q'-1\}$, such that
$$0 <\sigma(1).\rho-n_{\sigma(1)} < \sigma(2).\rho-n_{\sigma(2)}< \dots
<\sigma(q+q'-1).\rho-n_{\sigma(q+q'-1)} < 1;$$
\item  the integers 
$n_1,\dots,n_{q+q'-1}$ and the permutation $\sigma$ actually
do not depend on the choice of the number $\rho$ in 
$]\frac{p}{q}\,,\,\frac{p'}{q'}[$.
\end{itemize}
% Let $\rho$ be any number in the Farey interval 
%$]\,\frac{p}{q}\,,\,\frac{p'}{q'}\,[$. For
%every $k\in\ZZ$, we consider the integer $n_k$ such that 
%$k.\rho-n_k\in ]0,1[$. Then
%we consider the permutation $\sigma$ of the set $\{1,\dots,q+q'-1\}$, such that
%$$0 <\sigma(1).\rho-n_{\sigma(1)} < \sigma(2).\rho-n_{\sigma(2)}< \dots
%<\sigma(q+q'-1).\rho-n_{\sigma(q+q'-1)} < 1.$$
%It is relevant to notice here that the integers 
%$n_1,\dots,n_{q+q'-1}$ and the permutation $\sigma$ actually
%do not depend on the choice of the number $\rho$ in 
%$]\frac{p}{q}\,,\,\frac{p'}{q'}[$: this follows from
%proposition \ref{p.order} of appendix~\ref{app.arithmetic}.
  Then, for every
$k\in\{1,\dots,q+q'-1\}$, we  
consider the  homeomorphism
$\Phi_k:=\wt h ^{\sigma(k)}\circ T^{-n_k}$. Moreover, we set 
$\Phi_0:=\mbox{Id}$ and $\Phi_{q+q'}:=T$.
Finally, for every $k\in\{1,\dots,q+q'\}$, we consider the homeomorphism
$\Psi_k:=\Phi_k\circ\Phi_{k-1}^{-1}$ (see figure \ref{fig.arcs}).
% For every 
%$k$, the homeomorphism $\Psi_k$ is the
%product of a power of $\wt h$ and a power of $T$; in particular,
It is clear that each  $\Psi_k$ 
commutes with the translation $T$, and that these homeomorphisms are
pairwise commuting.

%\bigskip

Let $k\in\{2,\dots,q+q'-1\}$. We have $\Psi_k={\wt h}^{\sigma(k)}\circ
T^{-n_{\sigma(k)}}\circ {\wt h}^{-\sigma(k-1)}\circ T^{n_{\sigma(k-1)}}$. 
Hence, according to lemma \ref{lemm.rotation-itere},
the rotation set of  the homeomorphism $\Psi_k$ is:
$$
\mbox{Rot}(\Psi_k) =
\left\{ \sigma(k).\rho-n_{\sigma(k)}-\left 
(\sigma(k-1).\rho-n_{\sigma(k-1)}\right) \;\;
\mid \;\; \rho \in \rot(\wt h)
\; \right\}.
$$
Since by assumption the set $\rot(\wt h)$ is included in  $]
\frac{p}{q} \, , \, \frac{p'}{q'}[$, 
the definition of $\sigma$ implies that $\mbox{Rot}(\Psi_k)$
 is included in $]0,+\infty[$.
Similarly, we have:
$$
\mbox{Rot}(\Psi_1)=\left\{\sigma(1).\rho-n_{\sigma(1)} \;\; \mid
\;\; \rho\in  \rot(\wt h)  \right\},
$$
$$
\mbox{Rot}(\Psi_{q+q'})=\left\{1-\left(\sigma(q+q'-1).\rho-n_{\sigma(q+q'-1)} \right)
\;\; \mid \;\;
\rho\in   \rot(\wt h) \right\},
$$
and these sets are also included 
in $]0,+\infty[$.
Consequently, the homeomorphisms $\Psi_1,\dots,\Psi_{q+q'}$ satisfy the 
hypotheses of proposition~\ref{p.plusieurs-applications}. 
So this proposition provides us with an 
essential  simple arc
$\Gamma$ in the band $\widetilde\AA$ such that $R(\Gamma)$ is a 
strict attractor for each of the
homeomorphisms $\Psi_1,\dots,\Psi_{q+q'}$.
Let $\gamma$ be the projection in the annulus $\AA$ of the arc
$\Gamma$. It is an essential arc in the
annulus $\AA$~; let us prove that it is simple.
 To see this, we observe that the translation $T$ is equal
to the telescopic product $\Psi_{q+q'}\circ\dots\circ\Psi_1$. As a consequence, 
$R(\Gamma)$ is
a strict attractor for $T$. In
particular, the arc $\Gamma$ is disjoint from its image 
 $T$, and $\gamma$ is a simple arc.

\begin{figure}[ht]
\par \centerline{\hbox{\input{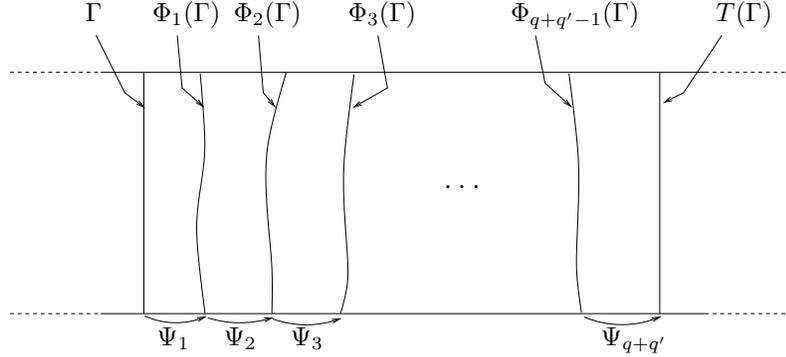}}} \par
\caption{\label{fig.arcs} The arcs
$\Gamma,\Phi_1(\Gamma),\dots,\Phi_{q+q'-1}(\Gamma)$ and $T(\Gamma)$}
\end{figure}

We are left to prove that the arcs
$\gamma,h(\gamma),\dots,h^{q+q'-1}(\gamma)$ are pairwise disjoint.
By construction of the arc $\Gamma$, we know that the set $R(\Gamma)$
is a strict attractor for the
homeomorphism $\Phi_1=\Psi_1$; in other words, the arc
$\Phi_1(\Gamma)$ is strictly on the right of the arc $\Gamma$.
For every $k\in\{2,\dots,q+q'-1\}$, we know that the set $R(\Gamma)$
is a strict attractor for the homeomorphism
$\Psi_k$;  since $\Psi_k$ and $\Phi_{k-1}$ commute, this implies that
the set $\Phi_{k-1}(R(\Gamma))$  is a strict
attractor for the homeomorphism $\Psi_k$; in other words, the arc
$\Psi_k\circ\Phi_{k-1}(\Gamma))=\Phi_k(\Gamma)$ is
strictly on  the right of the arc $\Phi_{k-1}(\Gamma)$.  Similarly,
%the set $\Phi_{q+q'-1}(R(\Gamma))$  is a strict
%attractor for the homeomorphism $\Psi_{q+q'}$; in other  words, 
the arc $\Psi_{q+q'}\circ\Phi_{q+q'-1}(\Gamma))=T(\Gamma)$  is strictly
on the right  of the arc $\Phi_{q+q'-1}(\Gamma)$
(see figure \ref{fig.arcs}). So we have proved that the arcs
$\Phi_1(\Gamma),\dots,\Phi_{q+q'-1}(\Gamma)$
are pairwise disjoint, and that all these arcs are strictly on the right of
$\Gamma$ and strictly on the left of $T(\Gamma)$,
\emph{i.e.} are included in the set $D:=R(\Gamma)\setminus
R(T(\Gamma))$. Since the set $D$ is a fundamental domain
for the covering map $\pi:\widetilde\AA\rightarrow\AA$, this implies that
the projections of the
arcs $\Gamma,\Phi_1(\Gamma),\dots,\Phi_{q+q'-1}(\Gamma)$ are pairwise
disjoint in the annulus $\AA$.
Now we observe that, for every $k$, the homeomorphism $\Phi_k$ is,
by definition, a lift of the homeomorphism
$h^{\sigma(k)}$; in particular, the projection of the arc
$\Phi_k(\Gamma)$ is the arc $h^{\sigma(k)}(\gamma)$. Hence,
we have proved that the arcs
$\gamma,h^{\sigma(1)}(\gamma),\dots,h^{\sigma(q+q'-1)}(\gamma)$ are
pairwise disjoint.
Since $\sigma$ is a permutation, this is equivalent to the fact that
the arcs $\gamma,h(\gamma),\dots,h^{q+q'-1}(\gamma)$
are pairwise disjoint.
\end{proof}

\subsection{Proof of the perturbation lemma}
\label{ss.perturbation}

\begin{proof}[Proof of lemma \ref{l.perturbation}] 
According to lemma \ref{lemm.ca-pousse}, there exists a positive integer $N$ such 
that the arc
$\Psi^{N}(\Gamma)$ is included in the interior of $R(\Gamma)$, 
\emph{i.e.} such that
$R(\Gamma)$ is a strict attractor for the homeomorphism $\Psi^{N}$. 
If $N=1$, then we
can set $\widehat\Gamma:=\Gamma$. Hence lemma \ref{l.perturbation} 
follows of sublemma \ref{sl.perturbation} below by induction on $N$.
\end{proof}

\begin{sublemma}
\label{sl.perturbation}
Let $\Psi:\widetilde\AA\rightarrow\widetilde\AA$ be a homeomorphism 
isotopic to the identity.
Suppose that we have found an essential simple arc $\Gamma$ in 
$\widetilde\AA$ and an integer
$n\geq 2$, such that $R(\Gamma)$ is an attractor for $\Psi$ and a 
strict attractor for $\Psi^n$.
Then, arbitrarily close to the arc $\Gamma$, there exists an 
essential simple arc $\Gamma'$ such
that $R(\Gamma')$ is an attractor for $\Psi$ and a strict attractor 
for $\Psi^{n-1}$.
\end{sublemma}

\begin{figure}[ht]
\par \centerline{\hbox{\input{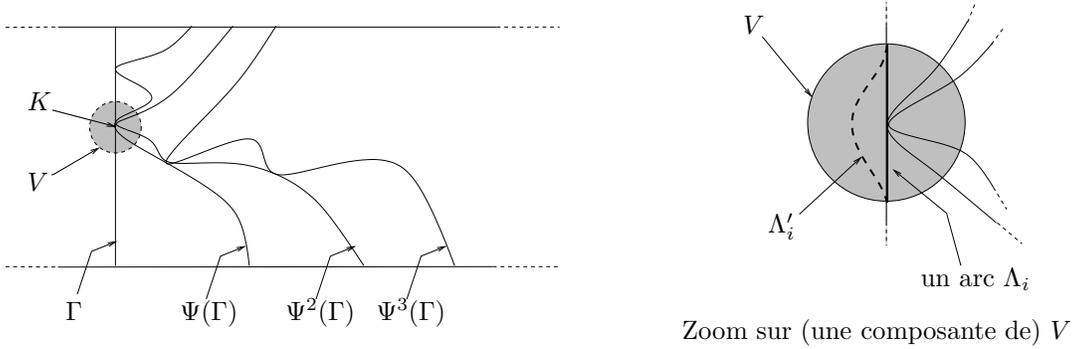}}} \par
\caption{\label{fig.gamma-prime}Construction of the arc $\Gamma'$ 
in sublemma \ref{sl.perturbation}}
\end{figure}

\begin{proof}
By Schoenflies theorem (\cite{Cai}), there exists a homeomorphism $G$ of the band 
$\widetilde\AA$, isotopic to the
identity, which maps the arc $\Gamma$ on the vertical segment 
$\{0\}\times [0,1]$.  So, up to replacing
$\Gamma$ by $G(\Gamma)$ and $\Psi$ by $G\circ\Psi\circ G^{-1}$,  we 
may assume that
$\Gamma$ is  the vertical segment $\{0\}\times [0,1]$.

Let us consider the compact set  $K:=\Psi^{n-1}(\Gamma)\cap\Gamma$. We have
$\Psi(K)\subset\Psi^n(\Gamma)\subset \inte(R(\Gamma))$ (since 
$R(\Gamma)$ is a strict attractor for
$\Psi^n$). Hence we can find a neighbourhood $V$ of $K$ such that 
$\Psi(V)\subset \inte(R(\Gamma))$.
Moreover, since $\Gamma$ is the vertical segment $\{0\}\times [0,1]$,
we can choose the neighbourhood $V$ such that 
$V\cap\Gamma$ is made of a finite number of
subarcs $\Lambda_1,\dots,\Lambda_p$ of the arc $\Gamma$. We construct 
an essential simple arc
$\Gamma'$ as follows: starting with the arc $\Gamma$, we replace each subarc 
$\Lambda_i$ by an arc $\Lambda_i'$,
which has the same ends as $\Lambda_i$, which is included in $V$, and 
whose interior is disjoint from
$R(\Gamma)$ (see figure \ref{fig.gamma-prime}). Observe that, by construction, we have
$\Gamma'\subset\Gamma\cup V$, $R(\Gamma)\subset R(\Gamma')$, 
$K\subset R(\Gamma')$, and $\Gamma'$ is arbitrarily close to $\Gamma$.

We have to prove first that the set $R(\Gamma')$ is an attractor for 
$\Psi$, \emph{i.e.} that the
arc $\Psi(\Gamma')$ is included in 
%%%the interior of
 $R(\Gamma')$. For 
that purpose, we recall that the arc
$\Gamma'$ is included in $\Gamma\cup V$,
% (by construction of $\Gamma'$),
 that the arc $\Psi(\Gamma)$ is
included in $R(\Gamma)$,
%(since $R(\Gamma)$ is an attractor for $\Psi$),
 that $V$ was chosen in such a
way that $\Psi(V)$ is included in the interior of $R(\Gamma)$, and 
that $R(\Gamma)$ is
included in $R(\Gamma')$. Hence we have
$$
\Psi(\Gamma')\subset\Psi(\Gamma)\cup \Psi(V)\subset
R(\Gamma)\subset R(\Gamma')
$$
\emph{i.e.} the set $R(\Gamma')$ is an attractor for $\Psi$.

Next we have to prove that the set $R(\Gamma')$ is a strict attractor 
for $\Psi^{n-1}$.
% \emph{i.e.} that the
%arc $\Psi^{n-1}(\Gamma')$ is included in the interior of 
%$R(\Gamma')$.
 For that purpose, we first recall that
%the arc $\Gamma'$ is included in $\Gamma\cup V$, and that
 $\Psi(V)$  is included in the
interior of $R(\Gamma)$. Since $R(\Gamma)$ is an attractor for 
$\Psi$, this implies that the set
$\Psi^{n-1}(V)$ is also included in the interior of $R(\Gamma)$. Now, we write
$\Psi^{n-1}(\Gamma)=\left(\Psi^{n-1}(\Gamma)\setminus K\right)\cup 
K$. On the one hand, we have
$\Psi^{n-1}(\Gamma)\setminus K\subset \Psi^{n-1}(\Gamma)\setminus 
\Gamma \subset
R(\Gamma)\setminus\Gamma=\inte(R(\Gamma))$ (the first inclusion 
follows from the definition of $K$, and
the second inclusion follows from the fact that $R(\Gamma)$ is an 
attractor for $\Psi$). On the other hand,
we recall that the set $K$ is included in the interior of 
$R(\Gamma')$ (by construction of the arc
$\Gamma'$). Hence, we have
$$\begin{array}{lllll}
\Psi^{n-1}(\Gamma') & \subset & \Psi^{n-1}(V)\cup\Psi^{n-1}(\Gamma)&&\\
 & \subset & \Psi^{n-1}(V)\cup\left(\Psi^{n-1}(\Gamma)\setminus
K\right)\cup K&&\\
  & \subset & \inte(R(\Gamma))\cup \inte(R(\Gamma)) \cup\inte(R(\Gamma'))
  & \subset & \inte(R(\Gamma')).
\end{array}$$
In other words, the set $R(\Gamma')$ is a strict attractor for 
$\Psi^{n-1}$. This completes the proof of sublemma~\ref{sl.perturbation}.
\end{proof}

%%%%%%%
\subsection{Proof of lemma \ref{l.vee}}
\label{s.lvee}
The proof  of lemma \ref{l.vee} relies heavily on the following classical
result of Ker\'ej\'art\'o (see \cite{Ker}, or \cite[page
246]{LeCYoc} for a modern proof): \emph{let $U_1$ and $U_2$ be two Jordan 
domains in the two-sphere, that is,  connected open sets whose boundary is 
homeomorphic to the circle~; then
each connected component of $U_1\cap U_2$ is also a Jordan domain}.

\begin{proof}[Proof of lemma \ref{l.vee}]
%So, we consider two essential simple arcs $\Gamma_1$ and $\Gamma_2$ 
%in $\widetilde\AA$, we
%denote by $E$ the unique non-compact connected component of the set 
%$(\widetilde\AA\setminus
%R(\Gamma_1))\cap (\widetilde\AA\setminus R(\Gamma_2))$,  and we 
Under the hypotheses of the lemma, let 
 $\Gamma$ denote the boundary of
$U$ in $\widetilde\AA$. We have to prove that $\Gamma$ is an 
essential simple arc. For that purpose, we
see the band $\widetilde\AA=\RR\times [0,1]$ as a subset $\RR^2$, and 
we see the
two-sphere $\SS^2=\RR^2\cup\{\infty\}$ as the one-point 
compactification of $\RR^2$. Then
Ker\'ej\'art\'o's result stated above implies that $U$ is a 
Jordan domain in $\SS^2$. In
particular, $\Gamma$ is included in a Jordan curve of $\SS^2$ (passing 
through the point $\infty$). From this it follows easily that $\Gamma$
 is an essential simple arc.
%Moreover, $\Gamma$ is a compact subset of $\RR^2$ (since $\Gamma$ is included
%in $\Gamma_1\cup\Gamma_2$). Hence, $\Gamma$ is a union of disjoined 
%arcs, the ends of these arcs being
%on the boundary of $\widetilde\AA$ (fact 1). Now, we observe that, 
%for $M$ large enough, we have
%$]-\infty,M]\times [0,1] \subset E \subset [M,+\infty[\times [0,1]$. 
%This implies that $\Gamma$
%contains at least one essential arc (fact 2).  Conversely, since $E$ 
%is a subset of $\widetilde\AA$ which is
%connected and non-compact, $\Gamma$ cannot contain two essential arcs 
%(fact 3). Finally, since $\Gamma$
%is included in $\Gamma_1\cup\Gamma_2$, the intersection of $\Gamma$ 
%with each of the
%two boundary component of $\widetilde\AA$ contains at most two points 
%(fact 4). Facts 1, 2, 3 and 4
%imply that $\Gamma$ is an essential simple arc.
\end{proof}

%intro : COMMUTING map
%%%%%%%%
\section{Alternative proof}
\label{s.alternative-proof}
This section is devoted to a second proof of the ``arc translation 
theorem" \ref{th.closed-annulus}.
It can be considered as a variation on  the first proof given in the 
previous section. We use two
independent arguments. The first one is a purely arithmetic argument, 
and tells that it is enough to find
an essential simple arc which is disjoint from its images under the two
``first-return maps'' $\Phi_1=T^{-p}\circ \wt h ^{q}$ and
$\Phi_2= T^{p'}\circ \wt h^{-q'} $. The second argument goes the 
following way. Suppose we are given a family of $k$ 
pairwise commuting 
maps, and consider sequences obtained by  starting with any point in the
closed band $\tA$ and iterating each time by one of the maps of the 
family (that is, we are
considering a positive orbit of the $\ZZ^k$-action generated by the
family). We prove that if the rotation sets of the $k$ maps are all
positive, then all the sequences obtained this way have a universally
bounded leftward displacement (actually, the proof is given only for
$k=2$, since we have the arithmetic argument in mind). Moreover, by 
continuity, this
remains true if we consider pseudo-orbits, \textit{i. e.} if a
little ``jump'' (or ``error'')  takes place at each step. Then we
construct the essential simple arc $\Gamma$ using a \emph{brick decomposition}.
This is a sort of triangulation which produces attractors in an
automatic way, as far as the behaviour of pseudo-orbits is controlled.
Brick decompositions were introduced by Flucher (\cite{Flu}). They have been used
by P. Le Calvez and A. Sauzet (\cite{LeCSau}) in order to prove the 
existence of attractors
for Brouwer homeomorphisms. In this text, we
only need the easy version, without the maximality property introduced
by A. Sauzet (\cite{Sau}).

%%%%%%%%
\subsection{Structure of the proof}
\label{ss.structure}
When $\Gamma_1$ and $\Gamma_2$ are two disjoint essential simple arcs 
in $\wt \AA$, there are two
possibilities: \\
-- either $\Gamma_2$ is ``on the right'' of $\Gamma_1$ and we write $\Gamma_1<\Gamma_2$,\\
-- or $\Gamma_1$  is ``on the right'' of $\Gamma_2$ and we write $\Gamma_2<\Gamma_1$.\\
In the appendix~\ref{app.arithmetic}, we shall prove the following result.

\begin{prop}
\label{p.arith}
Let $\wt h$ be the lift to $\wt \AA$ of some annulus homeomorphism 
$h$ which is isotopic to the identity
and let $\Gamma\subset \wt\AA$ be some essential simple arc in $\wt 
\AA$. Assume that there exists a
Farey interval $]\frac{p}{q},\frac{p'}{q'}  [$ such that
$$T^{-p'}\circ\wt h^{q'}(\Gamma)<\Gamma<T^{-p}\circ\wt h^{q}(\Gamma),$$
Then, the following properties hold:
\begin{enumerate}
\item  the arcs $T^{-\ell}\circ\wt h^{k}(\Gamma)$, with 
$k\in\{0,\dots,q+q'-1\}$ and $\ell\in\ZZ$, are
pairwise disjoint;
\item \label{e.disposition} these arcs are ordered in $\wt \AA$ as the lifts of the 
$q+q'-1$ first iterates of a vertical segment
in $\AA$ under the rotation $R_\alpha$ for any $\alpha\in 
]\frac{p}{q},\frac{p'}{q'}  [$. More precisely: given two
pairs of integers $(k,\ell)$ and $(k',\ell')$ in 
$\{0,\dots,q+q'-1\}\times\ZZ$, we have
$$
T^{-\ell}\circ\wt h^{k}(\Gamma)<T^{-\ell'}\circ\wt 
h^{k'}(\Gamma)\Longleftrightarrow
k \alpha - \ell < k' \alpha - \ell'.
%(k'-k)\alpha-(\ell'-\ell)>0.
$$
\end{enumerate}
\end{prop}

\begin{rema}\label{r.order}
In particular, the arc $\Gamma$ is disjoint from the arc 
$T^\ell(\Gamma)$ for any $\ell\in\ZZ$. Thus
$\gamma=\pi(\Gamma)$ is an essential simple arc in $\AA$, disjoint 
from its $q+q'-1$ first iterates
under $h$. Moreover, the cyclic order of the arcs 
$\gamma,h(\gamma),\dots,h^{q+q'-1}(\gamma)$ is the
same as the cyclic order of the iterates of a vertical segment under 
the rotation of angle $\alpha$, for any
$\alpha$ in the Farey interval $I=\left] \frac{p}{q},\frac{p'}{q'}\right[$.
\end{rema}

The above proposition will be combined with the following one.
\begin{prop}\label{prop.two-maps}
Let $\Phi_1$, $\Phi_2$ be two  homeomorphisms of
$\tA$, isotopic to the identity,  which
 commute and 
commute with the translation $T$, and whose rotation sets
  are included in $]0,+\infty[$.
Then there exists an  essential simple arc $\Gamma$ which is disjoint
from its images $\Phi_1(\Gamma)$ and $\Phi_2(\Gamma)$.
\end{prop}

Now we explain how theorem~\ref{th.closed-annulus}
 follows from these propositions. Then
the remaining of the section will be devoted to the proof of
proposition~\ref{prop.two-maps}.

\begin{proof}[Alternative proof of theorem~\ref{th.closed-annulus}
assuming propositions~\ref{p.arith} and~\ref{prop.two-maps}]
Let $\wt h$ be as in theorem~\ref{th.closed-annulus}. We consider the
two ``return maps'' $\Phi_1:=T^{-p}\circ \wt h ^{q}$ and $\Phi_2:= 
T^{p'}\circ \wt h^{-q'} $.
  According to lemma~\ref{lemm.rotation-itere}, the
rotation sets of both maps are included in $]0,+\infty[$. So we can
apply proposition~\ref{prop.two-maps}, and we get a curve $\Gamma$
which does not meet its images $\Phi_1(\Gamma)$ and $\Phi_2(\Gamma)$.
Then, since the rotation sets are positive, the order of
the curves must be such that  $\Phi_2^{-1}(\Gamma) < \Gamma <
\Phi_1(\Gamma)$ (this also follows from the proof of
proposition~\ref{prop.two-maps}). Now we can apply
proposition~\ref{p.arith}. Letting $\gamma$ be the projection of $\Gamma$ to
the annulus $\AA$, it follows that the arcs $\gamma$, $h(\gamma)$,
\dots , $h^{q+q'-1}(\gamma)$ are pairwise disjoint.
By remark~\ref{r.order} the cyclic order of these arcs is the same as 
the cyclic order of
the iterates of a vertical segment under the rigid rotation.
This concludes the proof of the theorem.
\end{proof}

%%%%%%%
\subsection{Pseudo-orbits for commuting homeomorphisms of
positive rotation sets}

During the whole section, we consider two homeomorphisms $\Phi_1$, $\Phi_2$ of
$\tA$ which 
commute and 
commute with the translation $T$, and we  make the assumption
that  both rotation sets of
$\Phi_1$ and $\Phi_2$ are included in $]0,+\infty[$.

A sequence $(x_n)_{n \geq 0}$
of points in $\tA$ is called a \emph{$(\Phi_1,\Phi_2)$-orbit} if
for all $n$, we have $x_{n+1}=\Phi_1(x_n)$ or $\Phi_2(x_n)$.
Let $d$ denote the Euclidean distance on $\tA = \RR \times [0,1]$ and 
$\varepsilon$
a positive real number.
  An \emph{$\epsilon$-$(\Phi_1,\Phi_2)$-pseudo-orbit} is a
sequence $(x_n)_{n \geq 0}$  of points in $\tA$ such that for all $n$,
$d(\Phi_1(x_n),x_{n+1}) < \epsilon$ or $d(\Phi_2(x_n),x_{n+1}) < \epsilon$.
The main result that makes this definition useful
is that  we can choose $\epsilon>0$
such that the leftward displacement of any 
$\epsilon$-$(\Phi_1,\Phi_2)$-pseudo-orbit is
universally bounded:
\begin{prop}
\label{prop.pseudo-orbits}
There exist $\epsilon >0$ and $M>0$ such that for any
$\epsilon$-$(\Phi_1,\Phi_2)$-pseudo-orbit $(x_n)_{n \geq 0}$,
for any $n \geq 0$,
$$p_1(x_n) \geq p_1(x_0)- M.$$
\end{prop}

To prove this proposition, we use lemma~\ref{l.deplacement} below 
which bounds  the leftward displacement of the 
$\varepsilon$-$(\Phi_1,\Phi_2)$-pseudo-orbits over long periods.
First, we prove a version of this lemma for $(\Phi_1,\Phi_2)$-orbits:
\begin{lemma}
\label{lemm.phi-orbites}
There exists an integer $N>0$ with the following property. For each
couple of non-negative integers $(N_1,N_2)$ such that $N_1+N_2 \geq 
N$, for every point
$x$ in $\tA$,
$$
p_1(\Phi_1^{N_1} \Phi_2^{N_2}(x)) \geq p_1(x)+2.
$$
\end{lemma}

\begin{proof}[Proof of lemma~\ref{lemm.phi-orbites}]
Applying lemma~\ref{lemm.ca-pousse} twice, we find numbers $\rho,s$
satisfying the inequality~(\ref{e.fact}) of this lemma for both $\Phi_1$ and
$\Phi_2$ (and for every point $x$ and every positive integer $n$).
Take a couple of non-negative integers $(N_1,N_2)$.
Then writing
$p_1(\Phi_1^{N_1} \Phi_2^{N_2}(x))- p_1(x)$ as the sum
$$
  \left[ p_1(\Phi_1^{N_1}
(\Phi_2^{N_2}(x))) - p_1(\Phi_2^{N_2}(x))\right] +
\left[ p_1(\Phi_2^{N_2}(x)) - p_1(x)  \right],$$
we see that this quantity is greater than
  $\rho (N_1+ N_2) -2 s$.
  We conclude that any integer $N$ such that $\rho N - 2 s \geq 2$
  will satisfy the conclusion of the lemma.
\end{proof}

Let us tackle the case of $\varepsilon$-$(\Phi_1,\Phi_2)$-pseudo-orbits:
\begin{lemma}\label{l.deplacement}
There exist an integer $N>0$ and a constant $\varepsilon>0$ with the 
following property. For every
$\varepsilon$-$(\Phi_1,\Phi_2)$-pseudo-orbit $(x_0,\dots,x_N)$ of length $N$,
$$
p_1(x_N) \geq p_1(x_0)+1.
$$
\end{lemma}

\begin{proof}[Proof of lemma~\ref{l.deplacement}]
Let $N$ be the integer given by lemma~\ref{lemm.phi-orbites}.
We shall say that an $\varepsilon$-$(\Phi_1,\Phi_2)$-pseudo-orbit 
$(x_0,\dots,x_N)$ of length $N$
is \emph{of type} $\sigma$, where $\sigma\in \{1,2\}^N$, if for every 
$n\in\{0,\dots,N-1\}$ we have
$d(x_{n+1},\Phi_{\sigma_{n+1}}(x_n))<\varepsilon$. Since the set 
$\{1,2\}^N$ is finite, it is sufficient to prove the lemma
for each type $\sigma$. In the remainder of the proof, the type 
$\sigma$ of the pseudo-orbits is fixed.

To prove the lemma we identify the tangent spaces $T_x \wt \AA$ with 
the plane $\RR^2$. Given a point
$x_0$ in $\wt \AA$ and a finite sequence of vectors of the plane
$\mathbf{v}=(\vec v_1,\dots, \vec v_{N})$, we define recursively
$$
x_1:=\Phi_{\sigma_1}(x_0)+\vec v_1\; , \; \dots\; , \; 
x_N:=\Phi_{\sigma_{N}}(x_{N-1})+\vec
v_{N}.
$$
 The vectors $\vec v_i$ will be chosen in the compact 
unitary ball $\DD$ of
$\RR^2$. Then we can consider the map
$$
\begin{array}{rcl}
{\cal F} : \tA \times (\DD)^N & \longrightarrow & \RR^2 \\
(x_0,\mathbf{v}) & \longmapsto  &  x_N.
\end{array}
$$
Clearly, the map $\cal F$ is continuous. Since ${\cal F}$ commutes
with the deck transformation $T$ (meaning that  ${\cal F}
(T(x_0),\mathbf{v})=T({\cal F} (x_0,\mathbf{v})$), and since the quotient
annulus $\AA$ is compact, $\cal F$ is uniformly continuous.
Therefore there exists $\epsilon \in ]0,1[$ such that for every $x_0$ in
$\tA$, for every sequence  $\mathbf{v}=(\vec v_1,
\dots, \vec v_{N})$ of vectors whose  Euclidean norms are less than $\epsilon$,
we have $d({\cal F} (x_0,\mathbf{v}),{\cal F} (x_0,(\vec 0))) < 1$.
We observe that~:
\begin{itemize}
\item  for every $\epsilon$-$(\Phi_1,\Phi_2)$-pseudo-orbit $(x_0,\dots,x_N)$ of
     type $\sigma$, $x_N$ can be expressed as ${\cal F}
     (x_0,\mathbf{v})$ for some  sequence  $\mathbf{v}=(\vec v_0,
\dots, \vec v_{N})$ with $\|\vec v_i \| < \epsilon$~;
\item we have the equality  ${\cal F} (x_0,(\vec 0))= 
\Phi_{\sigma_{N}}\circ \cdots \circ
\Phi_{\sigma_1} (x_0)$.
\end{itemize}
As a consequence, for every $\epsilon$-$(\Phi_1,\Phi_2)$-pseudo-orbit
$(x_0,\dots,x_N)$ of length
$N$ and type $\sigma$, we have the following inequalities:
$$p_1(x_N)-p_1(x_0) \geq p_1(\Phi_{\sigma_{N}}\circ \cdots \circ \Phi_{\sigma_1} 
(x_0))-p_1(x_0) -1 \geq 2 -1 =1.$$
(the latest inequality is a consequence of 
lemma~\ref{lemm.phi-orbites} applied to the
$(\Phi_1,\Phi_2)$-orbit $(x_0\, , \, \Phi_{\sigma_1}(x_0),\dots, 
\Phi_{\sigma_{N}}\circ \cdots \circ
\Phi_{\sigma_1} (x_0))$; 
here is the place where we use the fact that
$\Phi_1$ and $\Phi_2$ commute).
 This gives the lemma.
\end{proof}

We end with the proof of proposition~\ref{prop.pseudo-orbits}.

\begin{proof}[Proof of proposition~\ref{prop.pseudo-orbits}]
Let  $N$ and $\varepsilon>0$ be the integer and the constant given by 
lemma~\ref{l.deplacement}.
For any $\varepsilon$-$(\Phi_1,\Phi_2)$-pseudo-orbit 
$(x_0,\dots,x_\ell)$ of length $\ell<N$ we have,
\begin{equation}\label{e.casborne}
p_1(x_\ell)\geq p_1(x_0)-N(\varepsilon+s),
\end{equation}
where $s>0$ is the bound given by lemma~\ref{lemm.ca-pousse}.

Let us now consider any positive integer $n$ and any 
$\varepsilon$-$(\Phi_1,\Phi_2)$-pseudo-orbit $(x_0,\dots,x_n)$ of 
length $n$. We decompose $n$ as $kN+\ell$ with $k\geq 0$ and 
$\ell\in\{0\dots, N-1\}$. Lemma~\ref{l.deplacement} implies
$$p_1(x_{kN})\geq p_1(x_0)+k.$$
By~(\ref{e.casborne}) we get
$$p_1(x_{\ell})\geq p_1(x_{kN})-N(\varepsilon+s).$$
Putting all these inequalities together, one deduces
$$
\begin{array}{rcl}
p_1(x_n)-p_1(x_0)& = & [p_1(x_n) - p_1(x_{kN})] + [p_1(x_{kN})  - p_1(x_0)] \\
& \geq &  -N(\varepsilon+s).
\end{array}
$$
Hence, the proposition is proved for the constant $M=N(\varepsilon+s)$.
\end{proof}

\subsection{Brick decomposition}

\begin{figure}[ht]
\par \centerline{\hbox{\input{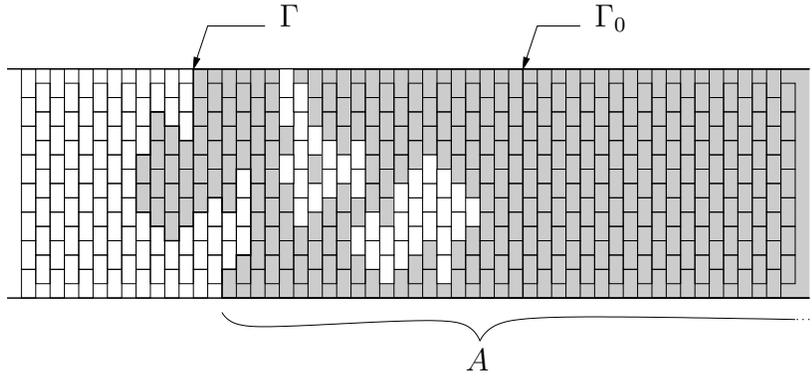}}} \par
\caption{\label{f.brick}A brick decomposition}
\end{figure}

We now turn to the proof of proposition~\ref{prop.two-maps}.
We consider a \emph{brick decomposition} of $\AA$, as shown on
figure~\ref{f.brick}.
%Essentially, this is a closed subset $F$ of $\tA$,
%containing the boundary of $\tA$~; it is locally homeomorphic to an
%arc, except at some points where it is homeomorphic to a
%\emph{triod}, \textit{i. e.} the union of three arcs having only one
%end-point in common.
Essentially, this amounts to taking an embedded triadic graph $F$ in $\tA$
  (\emph{triadic} meaning that each vertex belongs to exactly three
  edges). We demand that $F$ contains the boundary of $\tA$.
  A \emph{brick} is defined
to be the closure of a complementary domain of $F$ in $\AA$~; it is a
topological closed disk.
The last requirement in the definition of $F$ is the following key
  feature: every brick is of
diameter less than the number $\epsilon$ given by
proposition~\ref{prop.pseudo-orbits} (for the Euclidean metric on $\AA=\SS^1
\times [0,1]$).

\begin{rema}\label{r.manifold}
Note that, since $F$ is triadic,
the topological  boundary of the union of any family
of bricks is a 1-submanifold in $\tA$, with boundary included in the
boundary of $\AA$.
\end{rema}

A \emph{brick chain (from the brick $D_0$ to the brick $D_i$)} is a
sequence $(D_0, \dots D_i)$ of bricks in
$\tA$ such that $\Phi_1(D_0)\cup \Phi_2(D_0)$ meets $D_1$, \dots,
$\Phi_1(D_{i-1})\cup \Phi_2(D_{i-1})$  meets $D_i$.
%\begin{rema} The bricks are closed. Hence the interior of
%$\Phi_1(D_{k-1})\cup \Phi_2(D_{k-1})$  may be disjoint from the interior of
%$D_k$. We demand here that the closures intersect.
%\end{rema}

Take $\Gamma_0=\{0\} \times [0,1]$~; we can suppose that $\Gamma_0$ is
included in $F$ (as on figure \ref{f.brick}). We define a subset $A$
of $\tA$ in the following way: \\
-- to any brick $D_0$, we associate the union $\cD(D_0)$ of all the 
bricks $D$ of the
decomposition such that there exists a brick chain from $D_0$ to $D$;\\
-- the set $A$ is the union of all the sets $\cD(D_0)$, where $D_0$ 
ranges over the set
of all the bricks lying on the right of the arc $\Gamma_0$ (the brick 
$D_0$ may meet
$\Gamma_0$).

\begin{fact}\label{f.borneinf}
The set $A$ contains all the bricks on the
right of $\Gamma_0$ and is bounded to the left: there exists a 
constant $M$ such that
$A$ is included in $[-M,+\infty[\times [0,1]$.
\end{fact}
\begin{proof}
Indeed, if $(D_0, \dots D_i)$ is a brick chain, and $x$ is any point
  in $D_i$, then there exists an
  $\epsilon$-($\Phi_1$,$\Phi_2$)- pseudo-orbit $(x_n)$ such that $x_0$ is
  in $D_0$ and $x_i=x$. 
Remember that $\epsilon$ is given by
  proposition~\ref{prop.pseudo-orbits}; let $M$ be the other constant
  given by this proposition.
%Proposition~\ref{prop.pseudo-orbits} gives  $\varepsilon>0$
%that is used for the definition of the pseudo-orbits and a constant $M>0$.
Then we have  $p_1(x_i) \geq
  p_1(x_0)-M$, so that if $D_0$ is on the right of $\Gamma_0$,
  $p_1(x_0) \geq 0$, and
  $D_i$ is included in $[-M, +\infty [ \times [0,1]$. We conclude that
  $A$ is included  in $[-M, +\infty [ \times [0,1]$. The fact that $A$
contains all the bricks on the
right of $\Gamma_0$ follows from the definition of $A$ (considering
  chains made of  only one brick).
\end{proof}

\begin{fact}\label{fact.attractor}
The set $A$ is a strict attractor for $\Phi_1$ and $\Phi_2$, \textit{i. e.}
$$
\Phi_1(A) \subset \inte(A) \mbox{ and } \Phi_2(A) \subset \inte(A).
$$
\end{fact}
\begin{proof}
Indeed, let  $D$ be included in $A$. By definition, there exists a
brick chain   $(D_0, \dots D)$ with $D_0$ on the right of
$\Gamma_0$. Then for any brick $D'$ meeting $\Phi_1(D)$, the sequence
  $(D_0, \dots D, D')$ is again a brick chain, so $D'$ is also included
in $A$. Then the fact follows from the remark that $\Phi_1(D)$ is
included in the interior of the union of the bricks that it meets
(note that, to get a \emph{strict} attractor, it is crucial that the
bricks are defined to be \emph{closed}).
 Of course, the
same argument can be applied to the homeomorphism $\Phi_2$.
\end{proof}

Consider the essential arc $\Gamma$ ``bounding $A$ on the left'' (see 
figure \ref{f.brick}); more
precisely, using fact~\ref{f.borneinf} and remark~\ref{r.manifold}, 
this can be defined as the boundary of
the connected component of $\tA \setminus A$ containing
$]-\infty,-M[ \times [0,1]$. From fact~\ref{fact.attractor} it follows
that $\Gamma$ is disjoint from its images $\Phi_1(\Gamma)$ and
$\Phi_2(\Gamma)$. This ends the proof of proposition~\ref{prop.two-maps}.

%%%%%%%
%\subsection{End of the alternative proof}

%conjugacy classes
%\input{conjugacy-class-2}
\section{Closure of the conjugacy class of a pseudo-rotation: proof of corollary
\ref{c.closure-conjugacy-class}} 
\label{s.closure-conjugacy-class}

Using Poincar\'e's classical results (see the introduction), one can easily prove that, for any
orientation-preserving circle homeomorphism $h$ of rotation number $\alpha$, the rigid rotation of
angle $\alpha$ is in the closure of the conjugacy class of $h$. In
this section, we extend this result to 
irrational pseudo-rotations of the annulus, \emph{i.e.} we prove corollary
\ref{c.closure-conjugacy-class}. 

For every $\alpha\in\RR$, we denote by $T_\alpha$ the rigid translation in the
band $\wt \AA$ given by $ (\theta,t)  \mapsto
(\theta+\alpha,t)$. Of course, $T_\alpha$ is a lift of the rotation $R_\alpha$.

Corollary~\ref{c.closure-conjugacy-class} is an immediate consequence of the following proposition:

\begin{prop}\label{p.conjugaison} 
Let $h:\AA\rightarrow\AA$ be a homeomorphism that is isotopic to the identity and
$\wt h:\wt\AA\rightarrow\wt\AA$ be a lift of $h$. Suppose  that the
rotation set $\rot(\tilde h)$ is contained in
% the interior of
 some Farey interval $]\frac{p}{q},\frac{p'}{q'}[\subset \RR$. Then for any $\alpha\in
 ]\frac{p}{q},\frac{p'}{q'}[$,  there exists a
homeomorphism $\sigma$ of $\AA$, isotopic to the identity, such that for any lift $\wt \sigma$ of
$\sigma$ to $\wt \AA$, we have
$$\dd(\wt \sigma\circ \wt h\circ \wt \sigma^{-1}, T_\alpha)<\frac{30}{\min(q,q')}.$$
\end{prop}

\begin{proof}[Proof of corollary~\ref{c.closure-conjugacy-class} assuming
proposition \ref{p.conjugaison}] 
Let us consider an irrational pseudo-rotation 
$h$ of angle $\alpha$. Since $\alpha$ is irrational, it
belongs to some Farey interval $]\frac{p}{q},\frac{p'}{q'}[$ with $q$
and $q'$ arbitrarily large. Hence,
by proposition \ref{p.conjugaison}, $h$ is conjugate to some homeomorphism $\sigma\circ
h\circ \sigma^{-1}$ arbitrarily close to the rigid rotation $R_\alpha$.
\end{proof}

The remainder of section \ref{s.closure-conjugacy-class} is devoted to the proof of proposition
\ref{p.conjugaison}. Here is the idea of the proof. We begin by
applying our main theorem \ref{th.closed-annulus}, thus finding an
essential simple arc  $\gamma$ in $\AA$ which is disjoint from its first
$q+q'-1$ iterates under $h$. Let $\gamma_0$ be 
the vertical segment $\{0\} \times [0,1]$ in $\AA$
 (then $\gamma_0$ is disjoint from all its iterates
by the rotation $R_\alpha$).
 Since the cyclic order of the first iterates of $\gamma$ under $h$ is
the same as  the cyclic order of the first iterates of $\gamma_0$
under $R_\alpha$, one can perform a first conjugacy, by a homeomorphism
$\sigma_a$ sending $\gamma$ on $\gamma_0$,
 so that $h_a:=\sigma_a \circ h  \circ \sigma_a^{-1}$ coincides
with $R_\alpha$ on the iterates $\gamma_0, R_\alpha(\gamma_0), \dots,
R_\alpha^{q+q'-2}(\gamma_0)$ (first step of the proof). 

\begin{figure}[htp]
\par \centerline{\hbox{\input{farey2.pstex_t}}} \par
\caption{\label{f.farey}The dynamical tilings in $\AA$ and $\wt \AA$
(here $]\frac{p}{q},\frac{p'}{q'}[=]3/5,2/3[$)}
\end{figure}

In the second step, we use the dynamical tiling generated by these
arcs. More precisely, let us call $D$ and $D'$ the two tiles
 adjacent to the arc $\gamma_0$ in this tiling (see
figure~\ref{f.farey}). Since the $q'-1$ first iterates of $D$ under $h_a$ have
mutually disjoint interiors, and since they coincide with the iterates under
$R_\alpha$, we can conjugate $h_a$ by a homeomorphism  supported by
the union of these discs so that the conjugated homeomorphism
coincides with $R_\alpha$ on all 
these discs but the last one. We do the same on the iterates of
$D'$. Now we notice that the $q'-1$ first iterates of $D$, together
with the $q-1$ first iterates of $D'$, cover the whole annulus.
Thus the second step of the proof provides us with a
homeomorphism $\sigma_b$ so that $h_b:=\sigma_b \circ h_a  \circ
\sigma_b^{-1}$ coincides with $R_\alpha$ on the whole annulus except
on the set $R_\alpha^{q'-1}(D) \cup R_\alpha^{q-1}(D')$,
 which happens to be the topological disc $O':=R_\alpha^{-1}(D \cup
 D')$ (see figure~\ref{f.farey}).
Note that the interior of this disc $O'$ is disjoint from its first $s-1$ iterates,
where $s:=\min(q,q')$.

For the last step, we consider the difference homeomorphism
$g:=R_\alpha^{-1} \circ h_b$ on the topological disc $O'$. A key
lemma, dealing with 
disc homeomorphisms, allows us to write $g$  as
the composition of $N$
homeomorphisms $g_N, \cdots, g_1$ of the disc $O'$ which are $\epsilon$-close to the
identity, \emph{the integer $N$ depending on $\epsilon$ but not on
$g$}. We  choose $\epsilon$ so that $N<s$, hence the
disc $O'$ is disjoint from its first $N-1$ iterates, and we
can make a last conjugacy $\sigma_c$ that distributes the difference $g$ on these
iterates.  Thus we get a homeomorphism $h_c:=\sigma_c \circ h_b  \circ
\sigma_c^{-1}$ such that on $O'$ and its first
$N-1$ iterates, $h_c$ coincides with the rotation $R_\alpha$
up to one of the homeomorphisms $g_k$ 
(and consequently is $\epsilon$-close to $R_\alpha$),
 and such that $h_c$  still exactly  coincides with 
$R_\alpha$ everywhere else. Hence $h_c$  satisfies the conclusion of
proposition~\ref{p.conjugaison}.

%%%%%%%%%%%%%%%%%%%%%%%
\subsection{Preliminaries: decomposition of disc homeomorphisms}
\label{ss.disc-homeo}

We denote by $\DD$ the unitary closed disc for the Euclidean metric of $\RR^2$.
%$$\DD=B(0,1)=\{(x,y)\in\RR^2,x^2+y^2\leq 1\}.$$ 
We denote by $\homeo^+(\DD)$  the
set of homeomorphisms of the disc $\DD$ isotopic to the identity, and
by $\homeo(\DD,\partial\DD)$ the set of those that coincide
 with the identity on the boundary of $\DD$.
We consider the usual distance $d(h,h')=\sup\{d(h(x),h'(x)), x \in
\DD\}$ on these sets. We will say that two homeomorphisms $h,h'$ are
\emph{$\epsilon$-close} if $d(h,h') < \epsilon$.
The aim of this subsection is to prove the following lemma and corollary:

\begin{lemma}
\label{l.decomposition} 
For every $\varepsilon>0$, there exists $N\in \NN$ such that every homeomorphism
$h\in \homeo(\DD,\partial \DD)$ can be written as a product $h=h_N\circ\cdots\circ h_1$
of $N$ homeomorphisms 
%$h_1,\dots,h_N$ 
in $\homeo(\DD,\partial \DD)$ which are $\varepsilon$-close
to the identity.
Moreover, we can choose $N$  less than $\frac{4}{\epsilon}+4$.
\end{lemma}

Let $\alpha_1, \alpha_2$ be two disjoint closed arcs included in the boundary
$\partial \DD$. We denote by $\homeo(\DD,\alpha_1 \cup \alpha_2)$ the
set of homeomorphisms of the disc $\DD$ that coincide with the
identity on $\alpha_1 \cup \alpha_2$.
\begin{coro}
\label{c.decomposition}
For every $\varepsilon>0$, there exists $N\in \NN$ such that every homeomorphism
$h\in \homeo(\DD,\alpha_1 \cup \alpha_2)$ can be written as a product $h=h_N\circ\cdots\circ h_1$
of $N$ homeomorphisms 
%$h_1,\dots,h_N$ 
in $\homeo(\DD,\alpha_1 \cup
\alpha_2)$ which are $\varepsilon$-close to the identity.
Moreover, we can choose $N$  less than $\frac{6}{\epsilon}+5$.
\end{coro}

\begin{rema}
Note that, under the conclusion of the lemma, for every $k\leq N$, the homeomorphism $h_k\circ
h_{k-1}\circ\cdots\circ h_0$ is $\epsilon$-close to the homeomorphism
$h_{k-1}\circ\cdots\circ h_1$.
%: the distance between these two homeomorphism is less than $\varepsilon$.
 On the contrary, it is not true in general that the
homeomorphism $h_N\circ\cdots\circ h_k$ is $\epsilon$-close to the
homeomorphism $h_{N}\circ\cdots\circ h_{k+1}$.
\end{rema}

\begin{proof}[Proof of lemma~\ref{l.decomposition}]
We fix a number $\varepsilon>0$ and we consider a homeomorphism $h\in \homeo(\DD,\partial
\DD)$. 

\bigskip

\noindent \textit{Step 1. We use Alexander's trick to prove that $h$
can be written as a product $h=h'\circ h_0$ of two 
homeomorphisms $h_0,\,h'\in\homeo(\DD,\partial \DD)$, where $h_0$ is $\varepsilon$-close
to the identity, and where $h'$ coincides with the identity on a neighbourhood of $\partial\DD$.}

\smallskip

\noindent First, we extend $h$ on the whole plane $\RR^2$ by the identity on $\RR^2\setminus \DD$. 
Then for $t\in ]0,1]$, we consider the homeomorphism
$A_t\in\homeo(\DD,\partial \DD)$ defined by $A_t(x)=t.h(x/t)$ for every $x\in\DD$.
%By construction, we have $A_1=h$ and, for every $t\in [0,1[$, the homeomorphism $A_t$ coincides
%with the identity in a neighborhood of $\partial\DD$; furthermore, the
%map $t \mapsto A_t$ is continuous.
%Now, we pick some $t_0\in ]0,1[$ close to $1$,
%so that $d(A_{t_0}^{-1}\circ h)<\epsilon$, and we set $h_0:=A_{t_0}$ and $h':=A_{t_0}^-1\circ h$.
%The homeomorphisms $h'$ and $h_0$ satisfy all the desired properties.
Now we set $h':=A_{t_0}$ and $h_0:=A_{t_0}^{-1}\circ h$. It is easy to
check that if $t$ is close enough to $1$ the desired properties are satisfied.
\bigskip

\noindent \textit{Step 2. It remains to prove that the homeomorphism $h'$ can be
written as a product of $n$
elements of $\homeo(\DD,\partial\DD)$ which are $\epsilon$-close to
 the identity, with $n \leq \frac{4}{\varepsilon}+3$.}

\smallskip

\noindent Let $\delta>0$ be a number such that the homeomorphism $h'$ coincides with
the identity outside the Euclidean ball $B(0,1-\delta)$.
Let us consider now the radial homeomorphism $g_{\delta,\varepsilon}$
in $\homeo(\DD,\partial \DD)$ given by $$g_{\delta,\varepsilon}(x)=\varphi(\|x\|)x$$
where the homeomorphism $\varphi$ of $[0,1]$ is defined on figure~\ref{f.phi}.
\begin{figure}[htp]
\par \centerline{\hbox{\input{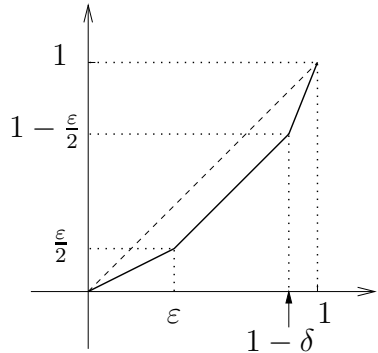}}} \par
\caption{\label{f.phi}The function $\varphi$}
\end{figure}
%$$\varphi(r)=
%\begin{cases}
%r-\frac{r}{2}, \quad \text{if $r\in[0,\varepsilon]$,}
%\\
%r-\frac\varepsilon 2, \quad \text{if $r\in [\varepsilon, 1-\delta]$,}\\
%r-\frac \varepsilon 2 \frac{1-r}{\delta},\quad \text{if $r\in [1-\delta,1]$.}
%\end{cases}
%$$
The homeomorphism $g_{\delta,\varepsilon}$ is $\varepsilon$-close to the identity.
Now, let $m$ be the integer in
$[\frac{2}{\varepsilon},\frac{2}{\varepsilon}+1[$. We observe that 
$\varphi^m(1-\delta)<\frac{\varepsilon}{2}$. Hence the homeomorphism
$g_{\delta,\varepsilon}^m$ maps  the ball $B(0,1-\delta)$ inside the ball $B(0,\varepsilon)$, 
and thus the homeomorphism $h_b=g_{\delta,\varepsilon}^m\circ h' \circ
g_{\delta,\varepsilon}^{-m}$ coincides with the identity  outside the
ball $B(0,\frac{\varepsilon}{2})$ (of center $0$ and radius $\frac{\varepsilon}{2}$).  
In particular, the homeomorphism $h_b$ is $\varepsilon$-close to the identity.
So, writing
$$h'=g_{\delta,\varepsilon}^{-m}\circ h_b \circ g_{\delta,\varepsilon}^{m}\;\;,$$
we obtain that $h'$ is the product of $n=2m+1$ elements of $\homeo(\DD,\partial\DD)$,
all of which are $\varepsilon$-close to the identity. Note that 
$n  \leq \frac{4}{\varepsilon}+3$ as announced.
\end{proof}

\begin{proof}[Proof of corollary~\ref{c.decomposition}]
We fix a number $\varepsilon>0$ and we consider a homeomorphism $h\in
\homeo(\DD,\alpha_1 \cup \alpha_2)$.

To begin with, note that every increasing homeomorphism of the interval $[0,1]$ can be
 written as a product of $N$ homeomorphisms that are $\epsilon$-close to the identity,
 with $N \leq \frac{1}{\epsilon}+1$.
From this we deduce that the same statement hold within the space 
$\homeo(\partial \DD,\alpha_1 \cup \alpha_2)$ of homeomorphism of the boundary circle
 $\partial \DD$ that are the identity on $\alpha_1 \cup \alpha_2$
(allowing $N \leq \frac{2}{\epsilon}+1$ because the diameter of this circle is $2$
since we use the metric on $\partial \DD$ as a subset of the disc $\DD$).
According to this we write
$$
h|_{\partial \DD}=H_N \circ \cdots \circ H_1.
$$
Now consider the ``circular extension mapping'' $\Phi$: using complex numbers notation,  
$$
\Phi : \left\{
\begin{array}{rcl}
\homeo(\partial \DD) & \longrightarrow & \homeo(\DD) \\
H & \longmapsto & \left( h : re^{i\theta}  \mapsto rH(e^{i\theta}) \right).
\end{array}
\right.
$$ 
Let $h':=\Phi(H_N)  \circ \cdots \circ \Phi(H_1)$, so that
$$
h=\Phi(H_N)  \circ \cdots \circ \Phi(H_1) \circ ({h'}^{-1} \circ h).
$$
In this formula each homeomorphism $\Phi(H_k)$ is $\epsilon$-close to the identity
 (since the circular extension mapping
  $\Phi$ is an isometry), and the homeomorphism ${h'}^{-1} \circ h$ is the identity 
on the boundary $\partial \DD$ (since $\Phi$ is an extension).
 Now we complete the proof by applying the previous 
lemma~\ref{l.decomposition} to the homeomorphism ${h'}^{-1} \circ h$.
\end{proof}

%%%%%%%%%%%%%%%%%%%%%%%%%%%%%%
%%%%%%%%%%%%%%%%%%%%%%%%%%%%%%
\subsection{Proof of proposition~\ref{p.conjugaison}}
\label{ss.proof-conjugaison}

Let us fix $\alpha$ in $]\frac{p}{q},\frac{p'}{q'}[$. We begin with some considerations on the rotation $R_\alpha$.

We consider the essential simple arcs $\Gamma_0=\{0\}\times [0,1]$ in $\wt\AA$, and 
$\gamma_0=\pi(\Gamma_0)$ in $\AA$. We consider the topological closed
discs (see figure~\ref{f.farey})
$$
\wt D:=R(\Gamma_0)\setminus \inte(R(T^{-p}\circ T_\alpha^q(\Gamma_0)))\quad
\mbox{ and }\quad
\wt D':=R(T^{-p'}\circ T_\alpha^{q'}(\Gamma_0))\setminus
\inte(R(\Gamma_0)).
$$ 
We set $D=\pi(\wt D)$ and $D'=\pi(\wt D')$, and we consider the family
of topological closed discs 
$$
\cD=\left\{R_\alpha^{k}(D)\;\;,\;\; k=0,\dots, q'-1\right\}\cup
\left\{R_\alpha^{k'}(D') \;\;,\;\; 
k'=0,\dots,q-1\right\}.
$$
\begin{lemma}\label{l.tiling}
The union of the discs in the family $\cD$ cover the annulus $\AA$, and their interiors are pairwise disjoint.
%The interiors of the elements of the family of discs $\cD$ are pairwise disjoint. The union of
%the elements of the family $\cD$ cover the annulus $\AA$.
\end{lemma}

\begin{proof}
By lemma \ref{p.not-contain-integer} of the appendix,
 no arc $R_\alpha^k(\gamma_0)$, with $k\in\{0,\dots,q+q'-1\}$, intersects 
the interior of $D$ nor $D'$. Consequently, the discs of the family $\cD$ have their interiors pairwise
disjoint.

Similarly, the arcs $R_\alpha^k(\gamma_0)$, with $k\in\{0,\dots,q+q'-1\}$, are pairwise disjoint. 
We observe that each of these arcs is contained in the boundary of exactly two elements of
the family $\cD$. This implies that  the family $\cD$ covers the annulus $\AA$.
\end{proof}

We denote by $O$ the topological closed disc $D'\cup D$. Observe that $O$ is also equal to
$R_\alpha^{q'}(D)\cup R_\alpha^q(D')$. Unless $q=q'=1$ (and in this case the proposition will
hold trivially from the step 1 below), $O$ is a topological closed disc.
The interior of $O$ is disjoint from its first $(\min(q,q')-1)$ backward iterates under $R_\alpha$.
The boundary of $O$ (as a topological manifold) is a simple closed curve $C$, and we have 
$$C=R_\alpha^q(\gamma_0)\cup R_\alpha^{q'}(\gamma_0)\cup C^+ \cup C^-,$$
where $C^-=\pi([q'\alpha-p',q\alpha-p]\times \{0\})$
and $C^+=\pi([q'\alpha-p',q\alpha-p]\times \{1\})$.

In order to compare the metrics on the topological disc $O$ and the
Euclidean disc $\DD$, we introduce a
homeomorphism $\psi : \DD\to O$: first, we note that $O$ is isometric to a Euclidean rectangle in
$\RR^2$, centered at $(0,0)$, with lengths $a,b\in ]0,1]$; then, we define $\psi$ as follows:
$$\psi(x,y)=\frac{\sqrt{x^2+y^2}}{\max(|x|,|y|)}.\left(\frac{ax}2,\frac{by}2\right).$$
A rough estimate shows that
\begin{equation}\label{e.contraction}
\forall z,z'\in\DD, \quad \dd_O(\psi(z),\psi(z'))\leq 2\dd_\DD(z,z'),
\end{equation}
where $\dd_O$ and $\dd_\DD$ are the usual Euclidean distances on $O$ and $\DD$.
Therefore one can easily transpose the result of
corollary~\ref{c.decomposition}, concerning the metric  space
$\homeo(\DD)$, to the space $\homeo(O)$ of homeomorphisms of $O$
isotopic to the identity endowed with the metric induced by the
Euclidean metric on $O$:
\begin{rema}\label{r.decomposition}
The statement in corollary~\ref{c.decomposition} still holds if we
replace the metric space  $\homeo(\DD,\alpha_1 \cup \alpha_2)$ by the
metric space $\homeo(O,R_\alpha^q(\gamma_0)\cup
R_\alpha^{q'}(\gamma_0))$ of homeomorphisms of $O$ which are the
identity on $R_\alpha^q(\gamma_0)\cup R_\alpha^{q'}(\gamma_0)$,
provided that we allow $N$ to be less than $2(\frac{6}{\epsilon}+5)$.
\end{rema}
Actually, we will cheat a little bit more by applying this remark to
$R_\alpha^{-1}(O)$ : note that this is all right  since $R_\alpha$ is an isometry.

%%%%%%%%%%%%%%%%%%%%%%%%%%%%%%%%%%%%%%%%
%%%%%%%%%%%%%%%%%%%%%%%%%%%%%%%%%%%%%%%
\subsubsection{First step} 

In this first step, we build a homeomorphism $\sigma_a:\AA\rightarrow\AA$,
isotopic to the identity, such that the homeomorphism $h_a=\sigma_a\circ h\circ\sigma^{-1}_a$ 
coincides with the rotation $R_\alpha$ on each
arc $R_\alpha^k(\gamma_0)$ with $k\in\{0,\dots,q+q'-2\}$.

Applying the arc translation
 theorem \ref{th.closed-annulus} to the homeomorphism $h$, we obtain an
essential simple arc $\gamma$ in $\AA$ such that the $q+q'-1$ first
 iterates of $\gamma$ under $h$ are 
pairwise disjoint. Now, we consider a homeomorphism $\sigma_a:\gamma\rightarrow\gamma_0$
which maps the endpoint $\gamma\cap\partial^+\AA$ (resp. $\gamma\cap\partial^- \AA$) on the
endpoint $\gamma_0\cap\partial^+ \AA$ (resp. $\gamma_0\cap\partial^+ \AA$). 
Since the first $q+q'-1$ iterates of $\gamma$ are pairwise disjoint, we may first 
extend $\sigma_a$ to the union of these iterates, deciding that on $h^k(\gamma)$, 
$\sigma_a:=R_\alpha^k\circ \sigma_\alpha|_{\gamma}\circ h^{-k}$. Now the key fact is that
according to proposition~\ref{p.arith}, 
the cyclic order of the iterates of $\gamma$ under $h$ is the same as for the iterates 
of $\gamma_0$ under $R_\alpha$. Consequently, thanks to a
repeated use of Schoenflies theorem (see for example \cite{Cai}), we may further extend $\sigma_a$ to a
homeomorphism of $\AA$ (isotopic to the identity).
For every $k\in\{0,\dots,q+q'-2\}$, the arc $h^k(\gamma)$ is
mapped by $\sigma_a$ on the arc $R_\alpha^k(\gamma_0)$, and the
 conjugate $h_a:=\sigma_a\circ h\circ
\sigma^{-1}_a$ of $h$ coincides with the rotation $R_\alpha$ on the arc $R_\alpha^k(\gamma_0)$.

We note that, since the homeomorphisms $\wt\sigma_a\circ \wt h\circ \wt\sigma^{-1}_a$ and
the translation $T_\alpha$ coincide on the arc $T^k(\Gamma_0)$ for every $k\in\ZZ$, we have
$$\dd(\wt \sigma_a\circ \wt h \circ \wt \sigma_a^{-1}, T_\alpha)\leq 1.$$
Hence, proposition \ref{p.conjugaison} already holds if $\min(q,q')\leq 30$. So, in the remainder of
the proof, we shall assume that $\min(q,q')$ is bigger than $30$.

%%%%%%%%%%%%%
\subsubsection{Second step} 

In this second step, we build a homeomorphism $\sigma_b\in \homeo(\mathbb{A})$
which is isotopic to the identity, and such that the conjugate $h_b:=\sigma_b\circ h_a\circ
\sigma^{-1}_b$ coincides with the rotation $R_\alpha$ everywhere except possibly on the 
topological disc $O':=R_\alpha^{-1}(O)=R_\alpha^{q'-1}(D) \cup R_\alpha^{q-1}(D')$.

According to the first step, for each $k \in \{0, \dots, q'-1\}$,  
we have $h_a^k(D)=R_\alpha^k(D)$; on this disc, we let $\sigma_b:=R_\alpha^k \circ h_a^{-k}$.
 We use the same formula on the disc $h_a^k(D')=R_\alpha^k(D')$ for  each $k \in \{0, \dots, q-1\}$.
Note that the intersection of any two such discs, if not empty, is an 
arc $R_\alpha^k(\gamma_0)$ with $k \in \{0, \dots, q+q'-1\}$, and that on these arcs the map
 $R_\alpha^k \circ h_a^{-k}$ is the identity. Thus these formulae are coherent, 
and define a homeomorphism $\sigma_b$ of $\AA$ isotopic to the identity. Now
 one easily checks that the homeomorphism $h_b$ defined as indicated above coincides with $R_\alpha$ 
on each disc $R_\alpha^k(D)$ with $k \in \{0, \dots, q'-2\}$ and   
$R_\alpha^k(D')$ with $k \in \{0, \dots, q-2\}$.
Using lemma~\ref{l.tiling}, we see that the union of all these discs 
cover the whole annulus but the set $O'$, and thus the second step is complete.

%%%%%%%%%%%%
\subsubsection{Last step} 
We denote by $s$ the minimum of $q$ and $q'$.
We consider the homeomorphism $g$ of the topological disc $O'$ defined by 
$g=R_\alpha^{-1} \circ h_b$. Note that $g$ is the identity on the boundary arcs 
$R_\alpha^{q-1}(\gamma_0)$ and $R_\alpha^{q'-1}(\gamma_0)$. Now we apply 
the results on decomposition of disc homeomorphisms: according to 
corollary~\ref{c.decomposition} and remark~\ref{r.decomposition}, we can write
$$
g=g_N \circ \cdots \circ g_1
$$
where each $g_i$ is a homeomorphism of $O'$ which is the identity on $R_\alpha^{q-1}(\gamma_0) 
\cup R_\alpha^{q'-1}(\gamma_0)$, and which is $\epsilon$-close to the identity of $O'$.
In applying this we choose $\epsilon=\frac{12}{s-12}$, so that
$$N \leq 2(\frac{6}{\epsilon}+5) \leq s-2$$
and one can check that $\epsilon \leq \frac{30}{s}$.

Now for each $k \in \{1, \dots, N\}$, we define $\sigma_c$ on the disc
$R_\alpha^{k}(O')$ by the formula
$$
\sigma_c:=R_\alpha^k \circ g_k \circ \cdots \circ g_1 \circ h_b^{-k}.
$$
We also let $\sigma_c$ be equal to the identity on the remaining of the annulus $\AA$.
It remains to check that these formulae are coherent, and that the homeomorphism
 $h_c:=\sigma_c\circ h_b\circ \sigma_c^{-1}$ is $\epsilon$-close to the identity.

Since $N \leq s-2$, the discs $R_\alpha(O'), \dots , R_\alpha^{N}(O')$ have their 
interiors pairwise disjoint; furthermore, the formulae defining $\sigma_c$ gives 
the identity on the boundary arcs $R_\alpha^\ell(\gamma_0)$
 of these discs. This proves that $\sigma_c$ is 
a well-defined homeomorphism of $O'$. 

Now observe first that 
$\sigma_c$ is equal to $R_\alpha^N\circ g \circ h_b^{-N}$ 
on the disc $R_\alpha^N(O')$. Since $N<s-1$, we have $h^N_b=R_\alpha^N\circ g$ on
$O'$ (according to the second step). This shows that $\sigma_c$ is the
identity on $R_\alpha^N(O')$.

For $k\in\{0,\dots,N-1\}$, one can check that the homeomorphism $h_c$ is
equal to $R_\alpha \circ \left(R_\alpha^{k}\circ g_{k+1} \circ R_\alpha^{-k}\right)$ 
on the disc $R_\alpha^k(O')$. Since $g_{k+1}$ is $\epsilon$-close to the identity, 
and since $R_\alpha$ is an isometry, we get that $h_c$ is $\epsilon$-close to $R_\alpha$ 
on this disc.
%Hence,
%using  (\ref{e.distance}), we obtain that, for every $k\in\{0,\dots,N-1\}$, the distance
%between the homeomorphism $h_c$ and the rotation $R_\alpha$ is bounded by $20/s$ on the disc
%$R_\alpha^k(O)$. 
Moreover, the homeomorphisms $h_b$ and $h_c$ coincide on
$E:=\AA\setminus\bigcup_{k=0}^{N-1}R_\alpha^k(O')$, since
$\sigma_c=\id$ on $R_\alpha^N(O')$. According to the second step, this proves that 
$h_c$ is equal to $R_\alpha$ on the set $E$. Thus the proof of the last step is complete.

\section{Proof of corollary \ref{c.without-periodic-orbit}}

As stated in the introduction, corollary 
\ref{c.without-periodic-orbit} relies on theorem
\ref{th.pseudo-rotation} and on a generalization of 
Poincar\'e-Birkhoff theorem obtained by
Bonatti  and Guillou. Let us first  recall the statement of 
Bonatti-Guillou's result:

\begin{theo}[\cite{Gui}]
\label{th.BonattiGuillou}
Let $h$ be a homeomorphism of the annulus $\AA$, which is isotopic to 
the identity.
Assume that $h$ does not have any fixed point. Then:

\smallskip

\noindent (i) either there exists a non-homotopically trivial simple 
closed curve in $\AA$, which is
disjoint from its image under $h$,

\smallskip

\noindent (ii) or there exists an essential simple arc in $\AA$ which 
is disjoint from its image under $h$.
\end{theo}

To prove corollary \ref{c.closure-conjugacy-class}, we also need two 
technical lemmas:

\begin{lemma}
\label{l.courbe-fermee}
Let $h$ be a homeomorphism of the annulus $\AA$, which is isotopic to 
the identity.
Assume that $h$ does not have any fixed point. Assume that there 
exists a positive integer $p$, and a
non-homotopically trivial simple closed curve $\sigma$  in $\AA$ 
which is disjoint from  its image
under $h^p$. Then there exists a non-homotopically trivial simple 
closed curve $\widehat\sigma$
in $\AA$ which is disjoint from its image under $h$.
\end{lemma}

\begin{proof}[Sketch of the proof]
We use the same kind of arguments as in the beginning of section 
\ref{s.first-proof}. For every
non-homotopically trivial simple closed curve $\sigma$ in $\AA$, we 
denote by $B(\sigma)$
the closure of the connected component of $\AA\setminus\sigma$ which 
is ``below $\sigma$'' (that is, which contains $\SS^1 \times \{0\}$).
Given two non-homotopically trivial simple closed curves $\sigma_1$ 
and $\sigma_2$ in $\AA$, the
boundary of the  connected component of 
$(\AA\setminus B(\sigma_1))\cap
(\AA\setminus B(\sigma_2))$ containing $\SS^1 \times \{1\}$
is a non-homotopically trivial simple 
closed curve, that we denote by
$\sigma_1\vee\sigma_2$ (the proof is the same as for lemma \ref{l.vee}).

Now, we consider the integer $p$ and the curve $\sigma$ given in the 
hypothesis of lemma
\ref{l.courbe-fermee}. We consider the non-homotopically trivial 
simple closed curve
$$\widehat\sigma:=\sigma\vee h(\sigma)\vee \dots\vee h^p(\sigma).$$
The same arguments as in the beginning of the proof of proposition 
\ref{p.plusieurs-applications} show that
$B(\widehat\sigma)$ is an attractor for $h$ (\emph{i.e.} 
the image of $B(\widehat\sigma)$
under $h$ is included $B(\widehat\sigma)$). Finally, using the same 
arguments as in the proof of lemma
\ref{l.perturbation}, we can perturb the curve $\widehat\sigma$ in 
such a way that $B(\widehat\sigma)$
becomes a strict attractor for $h$ (\emph{i.e.} in such a way that, 
after the perturbation, the image of
$B(\widehat\sigma)$ under $h$ is included in the interior of 
$B(\widehat\sigma)$). In particular, we
obtain a non-homotopically trivial simple closed curve 
$\widehat\sigma$ which is disjoint
from its image under $h$.
\end{proof}

\begin{lemma}
\label{l.rotation-non-entiere}
Let $h$ be a homeomorphism of the annulus $\AA$, which is isotopic to 
the identity. Assume that $h$
does not have any fixed point. Assume moreover that there exists an 
essential simple arc $\gamma$ in
$\AA$ which is disjoint from its image under $h$. Then, for every lift 
$\wt h$ of $h$, the rotation set of $\wt
h$ is disjoint from $\ZZ$.
\end{lemma}

\begin{proof}
Let us consider a lift $\wt h$ of $h$ and an essential simple arc $\Gamma$ in $\wt \AA$ which is a lift of $\gamma$. By changing the lift $\wt h$, it is sufficient to prove that $0$ is not contained in the rotation set of $\wt h$.

Our assumption on $\gamma$ implies that $\Gamma$ is disjoint from its image under $\wt h$. We will suppose for instance that $\Gamma\cap R(\wt h(\Gamma))=\emptyset$. This implies that $R(\Gamma)$ is a strict attractor for $\wt h$. The same is true for $R(T(\Gamma))$.

Recall that $\wt \AA$ is the band $\RR\times [0,1]$. We define a symmetry $\chi$ of $\RR^2$ by $\chi(x,y)=(x,-y)$. One can then extend $\wt h$ to an orientation-preserving homeomorphism $\wt h'$ of $\RR^2$ by setting for any $(x,y)\in \wt \AA$ and $n\in \ZZ$:
\begin{itemize}
\item
$\wt h'(x,-y)=\chi(\wt h(x,y))$;
\item
$\wt h'(x,y+2n)=\wt h'(x,y)+(0,2n)$.
\end{itemize}
Like $h$, the homeomorphisms $\wt h$ and $\wt h'$ do not have any fixed point.
Brouwer's theory on homeomorphisms of the plane
 asserts that any orbit of $\wt h'$ goes to infinity in $\RR^2$ (see \cite{Gui}).
Let us consider any point $z_0\in R(\Gamma)$.
The second coordinate of the orbit of $z_0$ by $\wt h$ remains in $[0,1]$. Since the orbit of $z_0$ is the same by $\wt h$ and by $\wt h'$, the modulus of the first coordinate along the positive orbit of $z_0$ by $\wt h$ must take arbitrarily large values. Using the fact that $R(\Gamma)$ is an attractor, we deduce that the first coordinate along the positive orbit of $z_0$ is bounded from below and is not bounded from  above. Hence, there exists $n(z_0)\geq 0$ such  that $\wt h^{n(z_0)}(z_0)\in R(T(\Gamma))$. Since $R(T(\Gamma))$ is an attractor, any iterate $\wt h^{n}(z_0)$ with $n\geq n(z_0)$ belongs
to $R(T(\Gamma))$.

Note that since $R(T(\Gamma))$ is a strict attractor, for any $z'$ close to $z_0$ and any $n>n(z_0)$,
the iterate $\wt h^n(z')$ belongs to $R(T(\Gamma))$. By compacity of the square $Q=R(\Gamma)\setminus \inte(R(T(\Gamma)))$, there exists some integer $n_0\geq 1$ such that for any point $z\in Q$ and any $n\geq n_0$, the iterate $\wt h^{n}(z)$ belongs to $R(T(\Gamma))$.

Since $\wt h$ and $T$ commute and since $R(\Gamma)\setminus R(T(\Gamma))\subset Q$ is a fundamental domain for the covering $\wt \AA\to \AA$, one deduces that the rotation set of $\wt h$ is included in
$[\frac{1}{n_0},+\infty[$ and does not contained $0$. This concludes the proof.
\end{proof}

Now, we are ready to prove corollary \ref{c.closure-conjugacy-class}:

\begin{proof}[Proof of corollary \ref{c.closure-conjugacy-class}]
Let $h$ be a homeomorphism as in the statement of corollary 
\ref{c.closure-conjugacy-class},
\emph{i.e.} a homeomorphism of the annulus $\AA$, which is isotopic 
to the identity, and which
does not have periodic point. Let $\wt h:\wt\AA\rightarrow\wt\AA$ be a 
lift of $h$.

\medskip

\noindent \emph{First case.} If there exists a non-homotopically 
trivial simple closed curve in $\AA$
which is disjoint from its image under $h$,  then we are done.

\medskip

\noindent \emph{Second case.} Now, we assume that there does not 
exist any non-homotopically
trivial simple closed curve in $\AA$ which is disjoint from its image under $h$.
In order to apply theorem \ref{th.pseudo-rotation}, we have to prove 
that $h$ is an irrational
pseudo-rotation. Let $p$ be an integer. Since $h$ does not have any 
periodic point, the homeomorphism
$h^p$ does not any fixed point. So we can apply theorem 
\ref{th.BonattiGuillou} to the homeomorphism
$h^p$. Moreover, lemma \ref{l.courbe-fermee} implies that there does 
not exist any
non-homotopically simple closed curve in $\AA$ which is disjoint from 
its image under
$h^p$. Hence, theorem \ref{th.BonattiGuillou} provides us with an 
essential simple closed arc in $\AA$
which is disjoint  from its image under $h^p$. Then, using lemma 
\ref{l.rotation-non-entiere}, we obtain
that the rotation set of the homeomorphism $\wt h^p$ is disjoint 
from $\ZZ$. Equivalently (see lemma
\ref{lemm.rotation-itere}), the rotation set of $\wt h$ is disjoint 
from $\frac{1}{p}.\ZZ$. So, we have
proved that, for every
$p\in\ZZ$, the rotation set of $\wt h$ is disjoint from 
$\frac{1}{p}.\ZZ$. Since the rotation set of $\wt h$
is an interval, this implies that it is a single irrational number; 
in other words, $h$ is an
irrational pseudo-rotation. Then, for every $n\in\NN$, theorem 
\ref{th.pseudo-rotation} provides us with
an essential simple arc $\gamma_n$, such that the arcs 
$\gamma_n,h(\gamma_n),\dots,h^n(\gamma_n)$ are
pairwise disjoint.
\end{proof}
\appendix
\section{Some elementary properties of Farey intervals}
\label{app.arithmetic}

\subsection{Farey intervals and rotations}
\label{ssapp.arithmetic}

Let us fix a Farey interval
$I=\left]\frac{p}{q},\frac{p'}{q'}\right[$. We choose a  number $\alpha$ in $I$.
  In this subsection, we prove that the cyclic order of the
$q+q'-1$ first iterates of any orbit under the circle rotation of
angle $\alpha$
does not depend on the choice of $\alpha$ in $I$
(proposition A.2 below).

\begin{lemma}
\label{p.not-contain-integer}
For every $k\in\{1,\dots,q+q'-1\}$, the interval
$\left]k.\frac{p}{q},k.\frac{p'}{q'}\right[$ does not contain any
integer.
\end{lemma}

\begin{proof}
Let $k$ be a positive integer such that the interval
$]kp/q\,,\,kp'/q'[$ contains an
%\footnote{Quelqu'un a rajoute (necessarily positive), mais je ne vois pas
%pourquoi ni a quoi ca sert ???}
 integer $\ell$.
Then the rational number $\ell/k$ is in the Farey interval
$]p/q,p'/q'[$.
Let us write
%in particular, we have
$p'/q'-p/q=(p'/q'-\ell/k)+(\ell/k-p/q)$; we have
$\ell/k-p/q=(\ell q-kp)/(kq)\geq 1/(kq)$ and
$p'/q'-\ell/k=(p'k-q'\ell)/(kq')\geq 1/(kq')$. 
Besides, since $]p/q,p'/q'[$ is a
Farey interval, we have $p'/q'-p/q=1/(qq')$.
Putting everything
together, we obtain $1/(kq)+1/(kq')\leq 1/(qq')$ which is equivalent
to $k\geq q+q'$. 
\end{proof}

Lemma \ref{p.not-contain-integer} implies that, for every
$k\in\{1,\dots,q+q'-1\}$, there exists a unique integer $n_k\in \ZZ$,
such that the interval $]\,k.\frac{p}{q}-n_k,\,k.\frac{p'}{q'}-n_k[$
is included in $]0,1[$. Given a number $\alpha$ in the
Farey interval $I$, and an integer $k\in\{1,\dots,q+q'-1\}$, we set
$\alpha_k:=k\alpha-n_k$. We can now prove the announced result:

\begin{prop}
\label{p.order}
The order of the numbers $\alpha_1,\dots,\alpha_{q+q'-1}$ in $]0,1[$
does not depend on the choice of $\alpha$ in the Farey interval $I$.
\end{prop}

\begin{proof}
Let $k_1$ and $k_2$ be two integers in $\{1,\dots,q+q'-1\}$, with
$k_1\neq k_2$. Then, using lemma~\ref{p.not-contain-integer},
we see that the difference $\alpha_{k_2}-\alpha_{k_1}$ is never null if $\alpha$
is in $I$. Since this quantity depends continuously on
$\alpha$, its sign does not depend on the choice of $\alpha$ in
$I$. This completes the proof.
\end{proof}

\subsection{Farey intervals and rational approximations}
\label{ss.reduction}

It is well-known that any Farey interval $I=\left] p/q,p'/q'\right[$
is associated with a finite sequence of rational numbers
$\left(\frac{p_n}{q_n}\right)_{1\leq n\leq n_0}$ which satisfies the
following properties:
\begin{itemize}
\item $\frac{p_1}{q_1}$ and $\frac{p_2}{q_2}=\frac{p_1}{q_1}+1$ are two consecutive integers;
\item For any $2\leq n\leq n_0-1$, there exists $a_{n+1}\in\NN\setminus\{0\}$ such that
\begin{equation}
\label{e.reduite}
p_{n+1}=a_{n+1}p_n+p_{n-1} \;\;\;\mbox{ and }\;\;\;
q_{n+1}=a_{n+1}q_n+q_{n-1};
\end{equation}
\item $\frac{p_{n_0-1}}{q_{n_0-1}}$ and $\frac{p_{n_0}}{q_{n_0}}$ are the two endpoints of the
Farey interval $I$, that is
$$\left]\frac{p}{q},\frac{p'}{q'}\right[=
\left]\frac{p_{n_0-1}}{q_{n_0-1}},\frac{p_{n_0}}{q_{n_0}}\right[
\mbox{ or } 
\left]\frac{p}{q},\frac{p'}{q'}\right[=
\left]\frac{p_{n_0}}{q_{n_0}},\frac{p_{n_0-1}}{q_{n_0-1}}\right[.$$
\end{itemize}
The rational numbers $\frac{p_1}{q_1},\dots,\frac{p_{n_0}}{q_{n_0}}$
are the \emph{common Farey approximations} of the elements of the
interval $I$.

\begin{rema}
It is convenient to add to the sequence $\{(p_n,q_n),1\leq n \leq
n_0\}$ a first term $(p_0,q_0)=(1,0)$ so that~(\ref{e.reduite}) holds
also with $n=1$ and $a_2=1$.
\end{rema}

\subsection{Iterates of essential simple arcs: proof of proposition~\ref{p.arith}}

In order to prove proposition \ref{p.arith}, the main task is to show
that all the iterates of the arc  $\Gamma$ involved in the statement are
disjoint from $\Gamma$. To prove this point, we will also need to collect some
 information about their order (but not all the information): this programme is
realized by lemmas \ref{l.farey1} and \ref{l.farey3}. Actually, the
whole result concerning the order will follow easily from the disjointness.

We use the relation ``$\Gamma_1<\Gamma_2$'' introduced at the
beginning of section~\ref{ss.structure}.
We assume that we are given a homeomorphism $h:\AA\rightarrow\AA$
isotopic to the identity, a lift $\wt h:\wt\AA\rightarrow\wt\AA$ of
$h$, an essential simple arc $\Gamma$ in $\wt\AA$, and a Farey
interval $I=\left] p/q,p'/q'\right[$ such that $T^{-p'}\circ\wt
h^{q'}(\Gamma)<\Gamma<T^{-p}\circ
\wt h^{q}(\Gamma)$. We have to prove that the arcs $T^{-\ell}\circ\wt h^{k}(\Gamma)$ for
$k\in\{0,\dots,q+q'-1\}$ and $\ell\in\ZZ$ are pairwise disjoint, and
ordered as the lifts of the first iterates of a vertical segment under the rotation
$R_\alpha$ for any $\alpha\in ]p/q,p'/q'[$.  The proof relies on some
arithmetical properties of Farey intervals, and on the fact that $T$
and $\wt h$ commute. Let us first introduce some notations:
\begin{itemize}

\item For every $(\ell,k)\in\ZZ^2$, we denote by $\Gamma(\ell,k)$ the curve $T^{-\ell}\circ \wt
h^k(\Gamma)$.

\item We denote by $\frac{p_1}{q_1},\dots,\frac{p_{n_0}}{q_{n_0}}$ (with $p_i\wedge q_i=1$) 
the common Farey approximations of the elements of the interval $I$
(see subsection
\ref{ss.reduction}). Moreover, we set $p_0:=1$ and $q_0:=0$.  

\item For $0\leq n\leq n_0$,  we denote by $\wt{h}_n$ the map $T^{-p_n}\circ \wt h^{q_n}$.

\item Given three essential simple arcs $\Gamma, \Gamma',\Gamma''$,  we shall say that 
\emph{$\Gamma$ separates $\Gamma'$ and $\Gamma''$}, if these three arcs are disjoint and satisfy
either $\Gamma'<\Gamma<\Gamma''$ or $\Gamma''<\Gamma<\Gamma'$.  
%This
%amounts to saying that these three arcs are disjoint and that
%$\Gamma$ intersects either the band
%$R(\Gamma')\setminus\inte(R(\Gamma''))$ or the band
%$R(\Gamma'')\setminus\inte(R(\Gamma'))$.
\end{itemize}

We will use intensively the fact that if $H$ is
a homeomorphism of $\wt \AA$ isotopic to the identity, then $H$
``preserves inequalities'':  if $\Gamma_1$
and $\Gamma_2$ are two disjoint essential simple arcs in $\wt \AA$,
the following  holds:
\begin{equation}
\label{e.inequality}
\Gamma_1 < \Gamma_2    \Longrightarrow      H(\Gamma_1) < H(\Gamma_2).
\end{equation}
This remark explains why the argument of the proof of the proposition
will be essentially one-dimensional.
% be formally identical to the
%proof of the corresponding statement for circle homeomorphisms.

\begin{lemma}
\label{l.farey1}
For any $n\in\{0,\dots,n_0-1\}$, the arc $\Gamma$ separates
$\Gamma(p_n,q_n)$ and $\Gamma(p_{n+1},q_{n+1})$. Furthermore, if $n
\neq 0$, then 
$\Gamma(p_{n+1},q_{n+1})$ separates $\Gamma$ and $\Gamma(p_{n-1},q_{n-1})$.
\end{lemma}

\begin{rema}\label{r.farey1}
The order of the whole family of arcs involved in the lemma is the following
(for simplicity, we assume that $n_0$ is even):
$$
T^{-1}(\Gamma)=\Gamma(p_0,q_0) < \Gamma(p_2,q_2) < \cdots <
\Gamma(p_{n_0},q_{n_0})
<  \Gamma 
$$
$$
\Gamma < \Gamma(p_{n_0-1},q_{n_0-1}) < \cdots < \Gamma(p_1,q_1) <
T(\Gamma).
$$
\end{rema}

\begin{proof}[Proof of lemma \ref{l.farey1}]
The proof is a decreasing induction on $n$. First, we observe that the
arc $\Gamma$ separates $\Gamma(p_{n_0-1},q_{n_0-1})$ and
$\Gamma(p_{n_0},q_{n_0})$ by assumption~(\ref{e.disposition})
 in proposition~\ref{p.arith}.

Now, suppose that we have proven that the arc $\Gamma$ separates the
arcs $\Gamma(p_{n},q_{n})$ and $\Gamma(p_{n+1},q_{n+1})$ for some
$n < n_0$. To fix ideas, we  assume that
\begin{equation}
\label{e.farey-1.1}
\wt{h}_n(\Gamma)<\Gamma<\wt{h}_{n+1}(\Gamma).
\end{equation}
The left ``inequality'' implies, using (\ref{e.inequality}) and the
 transitivity of ``$<$'',
 that $\Gamma < \wt{h}_n^{-a_{n+1}}(\Gamma)$.
Since $\wt{h}_{n+1}$ and $\wt{h}_n$ commute, using again (\ref{e.inequality}), we have 
\begin{equation}
\label{e.farey-1.2}
\wt{h}_{n+1}(\Gamma)<\wt{h}_{n}^{-a_{n+1}}\circ \wt{h}_{n+1}(\Gamma).
\end{equation}
By~(\ref{e.reduite}), we have $\wt{h}_{n-1}=\wt{h}_n^{-a_{n+1}}\circ \wt{h}_{n+1}$.
Hence, putting together (\ref{e.farey-1.1}) and (\ref{e.farey-1.2}),
we get the desired inequalities:
$$\Gamma(p_n,q_n)=\wt{h}_{n}(\Gamma) \; < \;  \Gamma \;  < \; 
 \wt{h}_{n+1}(\Gamma)=\Gamma(p_{n+1},q_{n+1})
 \; <  \; \wt{h}_{n-1}(\Gamma)=\Gamma(p_{n-1},q_{n-1}).$$
\end{proof}

\begin{proof}[Proof of remark \ref{r.farey1}]
The only inequality that is not contained in the lemma is the last
one. But we have $\Gamma(p_2,q_2) < \Gamma$, and composing with $T$
gives $T(\Gamma(p_2,q_2)) < \Gamma$; it remains to note that 
$T(\Gamma(p_2,q_2))=\Gamma(p_1,q_1)$ (see subsection~\ref{ss.reduction}).
\end{proof}

%Notice that for $n=0$, the lemma states that $\Gamma$ separates
%$T^{-1}(\Gamma)$ and $\Gamma(p_1,q_1)$, and in particular $\Gamma$ and
%$T^{-1}(\Gamma)$ are disjoint. This implies that the arcs $\Gamma$ and 
%$T^{-\ell}(\Gamma)$ are disjoint for every integer $\ell$. This will
%also be the case $n=0$ in the next lemma:
\begin{lemma}
\label{l.farey3}
For every integer $n\in\{0,\dots,n_0-1\}$, and every pair of integers
$(\ell,k)\neq (0,0)$ with $\ell\in\ZZ$ and $k\in
\{0,\dots,q_n+q_{n+1}-1\}$, the arcs $\Gamma$ and $\Gamma(\ell,k)$ are
disjoint. 
%Moreover, the arc $\Gamma(\ell,k)$ does not meet nor separate the arcs
%$\Gamma(p_n,q_n)$ and $\Gamma(p_{n+1},q_{n+1})$.
\end{lemma}

To avoid the multiplication of cases, we will use the notation
$\Gamma \leq \Gamma'$ to mean that $\Gamma' \subset R(\Gamma)$
(remember that $R(\Gamma)$ denotes 
the \emph{closure} of the connected component of $\widetilde\AA\setminus\Gamma$ 
which is ``on the right" of the arc
$\Gamma$).
\begin{proof}
We will actually prove the following statement :
\emph{For $n,\ell,k$ as in the lemma,  the arc $\Gamma(\ell,k)$ does
not meet the open topological disc
whose boundary in $\wt \AA$ is  
$\Gamma(p_n,q_n) \cup \Gamma(p_{n+1},q_{n+1})$.}
 This, together with
lemma \ref{l.farey1}, implies the lemma.
We proceed by induction on $n$.

 The case $n=0$ comes from the inequalities $T^{-1}(\Gamma) < \Gamma <
 h_1(\Gamma) < T(\Gamma)$ extracted from remark~\ref{r.farey1}
(note that $q_0+q_1-1=0$).

Now we assume that the statement holds for some integer $n-1\geq 0$. By
lemma~\ref{l.farey1}, we may assume for instance that
$\Gamma(p_{n},q_{n})<\Gamma<\Gamma(p_{n-1},q_{n-1})$.
Note that  this in turn implies
$\Gamma(p_{n},q_{n})<\Gamma<\Gamma(p_{n+1},q_{n+1})$ (using the same lemma).
 By the induction
hypothesis, every arc $\Gamma(\ell',k')$ with $(\ell',k')\neq (0,0)$,
$\ell'\in\ZZ$ and $k'\in
\{0,\dots,q_{n-1}+q_n-1\}$ satisfies either
$\Gamma(\ell',k') \leq \Gamma(p_{n},q_{n})$ or
$\Gamma(p_{n-1},q_{n-1}) \leq \Gamma(\ell',k')$. 

%does intersect the band $R(\Gamma(p_{n},q_{n}))\setminus
%\inte(R(\Gamma(p_{n-1},q_{n-1})))$.

Let us consider some arc $\Gamma(\ell,k)$ with $(\ell,k)\neq (0,0)$,
$\ell\in\ZZ$ and $k\in \{0,\dots,q_n+q_{n+1}-1\}$.  We have to prove
that either $\Gamma(\ell,k) \leq \Gamma(p_{n},q_{n})$ or
$\Gamma(p_{n+1},q_{n+1}) \leq \Gamma(\ell,k)$. 
According to lemma \ref{l.farey1}, we have $\Gamma(p_{n+1},q_{n+1})
\leq \Gamma(p_{n-1},q_{n-1})$; combining this with 
 our induction hypothesis solves the case $k \leq 
q_{n-1}+q_n-1$. 

We now suppose $k \geq 
q_{n-1}+q_n$. Then we write $k=sq_n+k'$ with
$s\in\{1,\dots,a_{n+1}\}$ and
$k'\in\{0,\dots,q_{n-1}+q_{n}-1\}$. Similarly, we write
$\ell=sp_n+\ell'$. By the induction hypothesis,
 one of the three following cases holds:
\begin{itemize}
\item $(\ell',k')=(0,0)$. In this case, we have 
$\Gamma(\ell,k)=\Gamma(sp_n,sq_n)=\wt{h}_n^s(\Gamma) \leq \wt{h}_n(\Gamma)=\Gamma(p_n,q_n)$
%Hence $\Gamma(\ell,k)<\Gamma(p_n,q_n)$.
%the arc $\Gamma(\ell,k)$ does not intersect $R(\Gamma(p_n,q_n))$.
(the inequality comes from the hypothesis $\wt{h}_n(\Gamma) <
\Gamma$ using (\ref{e.inequality}) and the transitivity of $<$).

\item $\Gamma(\ell',k')<\Gamma$. In this case
$\Gamma(\ell,k)=\wt{h}_n^s(\Gamma(\ell',k'))<\wt{h}_n^s(\Gamma)$
%T^{-\ell'}\circ \wt h^{k'}(\Gamma(p_n,q_n))<\Gamma(p_n,q_n)$.
 and we conclude using the previous case.

\item $\Gamma< \Gamma(\ell',k')$. In this case, using our induction
hypothesis, we obtain $\Gamma(p_{n-1},q_{n-1}) \leq \Gamma(\ell',k')$.
It follows that
$$
\begin{array}{rcl}
\Gamma(p_{n+1},q_{n+1}) = \wt{h}_{n-1} \circ
\wt{h}_n^{a_{n+1}}(\Gamma) & \underset{(a)}{\leq} & \wt{h}_{n-1} \circ
\wt{h}_n^{s}(\Gamma)\\
& = & \wt{h}_n^{s}(\Gamma(p_{n-1},q_{n-1}))  \\
&  \underset{(b)}{\leq} & \wt{h}_n^{s}(\Gamma(\ell',k'))=\Gamma(\ell,k)
\end{array}
$$
where we used for $(a)$ that  $\wt h_n(\Gamma) < \Gamma$
 and  $s \leq a_{n+1}$ and for $(b)$ the above inequality.
\end{itemize}
This completes the proof of the lemma.
\end{proof}

We now turn to the proof of proposition~\ref{p.arith}.
Let us first prove that the arcs are disjoint (first item). For this it
suffices to note that  
$$
\Gamma(\ell,k) \cap \Gamma(\ell',k')= T^{\ell'} \wt h
 ^{-k'} (\Gamma(\ell-\ell',k-k') \cap \Gamma) 
$$
 and to apply  the last lemma (for $n=n_0-1$). For
the second item we note that 
 $\Gamma(\ell,k) < \Gamma(\ell',k') \Leftrightarrow
\Gamma(\ell-\ell',k-k') < \Gamma$.
 It is now sufficient to show that, for any
$k\in\{0,\dots,q+q'-1\}$ and for any $\ell$,
\begin{equation*}
\begin{split}
T^{-\ell}\circ\wt h^k(\Gamma)<\Gamma \quad & \mbox{ implies that }
\quad k \frac{p'}{q'}-\ell \leq 0,\\
\Gamma < T^{-\ell}\circ\wt h^k(\Gamma) \quad & \mbox{ implies that }
\quad k \frac{p}{q}-\ell \geq 0.
\end{split}
\end{equation*}
Let us prove, for instance, the second implication. 
We suppose that $\Gamma < T^{-\ell}\circ\wt h^k(\Gamma)$.
First note that by assumption we have $T^{-p'}\circ\wt
h^{q'}(\Gamma)<\Gamma$, so that $\frac{p'}{q'} \neq \frac{l}{k}$.
This also implies that the rotation set of the map $T^{-p'}\circ\wt
h^{q'}$ is included in $]-\infty,0]$, so that the rotation set of $\wt
h$ is included in $]-\infty,\frac{p'}{q'}]$ by
lemma~\ref{lemm.rotation-itere}. Similarly, 
from the inequality $\Gamma<T^{-\ell}\circ\wt h^k(\Gamma)$ we get that
the rotation set of $\wt h$ is included in $[\frac{\ell}{k}, +\infty[$.
From all this we conclude that $k\frac{p'}{q'} - \ell \geq  0$.
 Since the integer $k$ is in $\{0,\dots,q+q'-1\}$, we get from
lemma~\ref{p.not-contain-integer} that $k\frac{p}{q}-\ell \geq 0$. This
completes the proof of proposition \ref{p.arith}.
 %app.arithmetic
%
%%%%%
\section{The case of the torus}
\label{app.torus}
In a way similar to section~\ref{ss.rotation}, the rotation set may be defined for a lift $\tilde h$ of a torus homeomorphism
that is isotopic to the identity. It is a convex compact set $\rot(\tilde h)$ of $\RR^2$.

We recall the statement of the theorem proved (for diffeomorphisms) by
J. Kwapisz~(\cite{Kwa1}) on the torus, and explain how our methods work in this
case.
\begin{theo}
\label{th.torus}
Let $h:\TT^2\rightarrow\TT^2$ be a homeomorphism of the torus
isotopic to the identity, and 
let $\widetilde h:\widetilde{\TT^2} \rightarrow\widetilde{\TT^2}$ be a lift of
$h$. Assume that the rotation set of
$\widetilde h$ is included  in a Farey band
$]\,\frac{p}{q},\frac{p'}{q'}\,[ \times \RR$. Then there exists an essential
simple closed curve $\gamma$ in
the torus  $\TT^2$ such that the curves
$\gamma,h(\gamma),\dots,h^{q+q'-1}(\gamma)$ are pairwise disjoint.
\end{theo}
Here an \emph{essential simple closed curve} is a topological circle
in $\TT^2$ homotopic to the circle $\{0\} \times \TT^1$, and a
\emph{Farey band} is a band $]\,\frac{p}{q},\frac{p'}{q'}\,[ \times
\RR$ where $]\,\frac{p}{q},\frac{p'}{q'}\,[$ is a Farey interval.

The proofs of the annulus case given in this paper work with very few
changes in the torus setting. However they do not take place on the universal covering
of the torus but on the intermediate covering $\RR\times \SS^1\to \SS^1\times \SS^1=\TT^2$.

The essential simple arcs are replaced by the notion of essential simple closed curve in $\TT^2$
or $\RR\times \SS^1$. Note also that the only information we need on the rotation set
is its first projection $p_1(\rot(\tilde f))$. This compact interval plays the role of the rotation set
in the annulus case and has to be compared to the Farey intervals of $\RR$.
 %app.torus
%
%\input{}
%
%\input{}


\begin{thebibliography}{150}
\bibitem{Cai}
Cairns, Stewart S.
An elementary proof of the Jordan-Schoenflies theorem.
\textit{Proc. Amer. Math. Soc.} \textbf{2} (1951), 860--867.

\bibitem{Den}
Denjoy, Arnaud.
Sur les courbes d\'efinies par les \'equations diff\'erentielles \`a la surface du tore.
\textit{J. Math. Pur. Appl.}, IX. Ser. \textbf{11} (1932), 333--375.

\bibitem{Flu}
Flucher, Martin.
Fixed points of measure preserving torus homeomorphisms.
\textit{Manuscripta Math.} \textbf{68}  (1990),  no. 3, 271--293.

\bibitem{Fra}
Franks, John.
Generalizations of the Poincar\'e-Birkhoff theorem. 
\textit{Ann. of Math.} (2)  \textbf{128}  (1988),  no. 1, 139--151.

\bibitem{Gui}
Guillou, Lucien.
Th\'eor\`eme de translation plane de Brouwer et g\'en\'eralisations du th\'eor\`eme de Poincar\'e-Birkhoff.
\textit{Topology} \textbf{33} (1994), no 2, 331--351.

\bibitem{Han}
Handel, Michael. 
A pathological area preserving $C\sp{\infty }$ diffeomorphism of the plane. 
\textit{Proc. Amer. Math. Soc.} \textbf{86} (1982), no. 1, 163--168.

\bibitem{Her}
Herman, Michael.
Construction of some curious diffeomorphisms of the Riemann sphere.
\textit{J. London Math. Soc.} (2)  \textbf{34}  (1986),  no. 2, 375--384.

\bibitem{Ker}
B\'ela Ker\'ekj\'art\'o.
\textit{Topologie (I).} 1923.

\bibitem{Kwa1}
Kwapisz, Jaroslaw.
A priori degeneracy of one-dimensional rotation sets for periodic point free torus maps.
\textit{Trans. Amer. Math. Soc.} \textbf{354} (2002), no 7, 2865--2895.

\bibitem{Kwa2}
Kwapisz, Jaroslaw.
Combinatorics of torus diffeomorphisms.
\textit{Ergodic Theory Dynam. Systems} \textbf{23} (2003), no.2, 559--586. 

\bibitem{LeCSau}
Le Calvez, Patrice and Sauzet, Alain.
Une d\'emonstration dynamique du th\'eor\`eme de translation de Brouwer.
\textit{Exposition. Math.} \textbf{14} (1996), no.3, 277--287.

\bibitem{LeCYoc}
Le Calvez, Patrice and Yoccoz, Jean-Christophe.
Un th\'eor\`eme d'indice pour les hom\'eomorphismes du plan au voisinage d'un point fixe.
\textit{Ann. of Math.} (2) \textbf{146} (1997),  no. 2, 241--293.

\bibitem{MisZie}
Misiurewicz, Michal and Ziemian, Krystyna.
Rotation sets for maps of tori.
\textit{J. London. Math Soc.} (2) \textbf{40} (1989), no.3, 490--506.

\bibitem{Poi}
Poincar\'e, Henri.
M\'emoire sur les courbes d\'efinies par une \'equation diff\'erentielle.
\textit{J. de Math Pures Appl.} S\'erie III \textbf{7} (1881), 375--422, \textbf{8} (1882), 251--296,
\textit{J. de Math Pures Appl.} S\'erie IV \textbf{1} (1885), 167--244, \textbf{2} (1886), 151--217.

 Also in \textit{{\OE}uvres de Henri Poincar\'e}, tome I, Gauthier Villars, Paris (1928), 3--44, 44--84,
90--158 and 167--222.


\bibitem{Sau}
Sauzet, Alain.
\textit{Application des d\'ecompositions libres \`a l'\'etude des hom\'eomorphismes de surfaces.}
PhD thesis, Universit\'e Paris Nord, 2001.

\bibitem{Sch}
Schwartzman, Solomon. 
Asymptotic cycles.
\textit{Ann. of Math.} (2) \textbf{66} (1957), 270--284.

\end{thebibliography}
\end{document}